\documentclass[11pt,letter]{article}

\usepackage{framed}
\usepackage{setspace}
\usepackage{fullpage}
\usepackage[top=2cm, bottom=2.5cm, left=2cm, right=2cm]{geometry}
\usepackage{tkz-berge}
\newcommand{\argmax}[1]{\underset{#1}{\text{argmax}}}
\newcommand{\all}{\forall}
\newcommand{\olu}{\overline{u}}
\newcommand{\olo}{\overline{o}}
\newcommand{\ulu}{\underline{u}}
\newcommand{\ulo}{\underline{o}}
\newcommand{\argmin}[1]{\underset{#1}{\text{argmin}}}
\newcommand{\x}[2]{x^{#1}_{#2}}

\usepackage{mathtools}

\usepackage{graphicx}
\usepackage{algorithm,lipsum,changepage}
\usepackage[noend]{algpseudocode}
\usepackage{fullwidth}
\usepackage{enumitem}
\usepackage{pifont}
\usepackage{verbatim}
\usepackage{amsmath}
\usepackage{amsthm}
\usepackage{amssymb}
\usepackage{fix-cm}
\begin{document}
	\title{Robust Appointment Scheduling with Heterogeneous Costs}
\author{ Andreas S. Schulz \\
	Technische Universit\"at M\"unchen\\ andreas.s.schulz@tum.de \and  \and Rajan Udwani \\ UC Berkeley \\ rudwani@berkeley.edu} 
	\date{}
\maketitle
%


%
%
\begin{abstract}
	Designing simple appointment systems that under uncertainty in service times, try to achieve both high utilization of expensive medical equipment and personnel as well as short waiting time for patients, has long been an interesting and challenging problem in health care. We consider a robust version of the appointment scheduling problem, introduced by Mittal et al.\ (2014), with the goal of finding simple and easy-to-use algorithms. Previous work focused on the special case where per-unit costs due to under-utilization of equipment/personnel are homogeneous i.e., costs are linear and identical. We consider the heterogeneous case and devise an LP that has a simple closed-form solution.  This solution yields the first constant-factor approximation for the problem. We also find special cases beyond homogeneous costs where the LP leads to closed form optimal schedules. Our approach and results extend more generally to convex piece-wise linear costs.
	
	For the case where the order of patients is changeable, we focus on linear costs and show that the problem is strongly NP-hard when the under-utilization costs are heterogeneous. For changeable order with homogeneous under-utilization costs, it was previously shown that an EPTAS exists. We instead find an extremely simple, ratio-based ordering that is 1.0604 approximate. 

\end{abstract}

\newtheorem{lem}{Lemma}
\newtheorem{coro}[lem]{Corollary}
\newtheorem{thm}[lem]{Theorem}
\newtheorem{cnj}[lem]{Conjecture}
\newtheorem{eg}{Example}[section]
	\newtheorem{fact}[lem]{Fact}
		\newtheorem{proposition}[lem]{Proposition}
\section{Introduction} \label{intro}
Consider the problem of scheduling appointments in service operations where customers are served sequentially by a single server. Service times of customers are uncertain, and we wish to assign time slots for serving the customers in advance. An important practical setting where this problem arises everyday is in health care services, where there are numerous instances that require efficient scheduling of appointments, such as scheduling outpatient appointments in primary care and specialty clinics, scheduling surgeries for operating rooms, or appointments for MRI scans. 	
Often in these settings the order in which patients undergo the procedure is known in advance. Then a day before the procedures, a hospital manager determines planned start times and how much time to allot to each procedure. If the manager allots too small an interval to a procedure, it could easily go overtime and delay the next procedure. The inconvenience and costs resulting from such a delay are referred to as the \emph{overage} cost for that procedure. On the contrary, if the manager assigns a very large interval, then the procedure will likely end much sooner. In this case the hospital incurs an \emph{underage} cost as the equipment and personnel may be left idle until the scheduled start of the next procedure. 
We would like to design an appointment schedule that can achieve the desired trade-off between overage and underage costs.

 In an influential survey on outpatient scheduling, Cayirli and Veral \cite{cayi}
	concluded that the ``biggest challenge for future research will be to develop easy-to-use heuristics''. Traditionally, most models for the problem are stochastic in nature and assume distributional knowledge of the uncertain service times. While such models are very powerful, estimating distributions accurately can require large amounts of data. There are many settings where such data is in fact available but, in applications related to health care appointment scheduling there is evidence (ref.\ \cite{macar,mak,denton}) that the amount of data available by surgery types, let alone by surgery types and surgeons, is extremely limited. Moreover, computing an objective function that involves finding the expectation of a high dimensional non-linear function can be computationally burdensome. Further, Mak et al.\ \cite{mak} and Mittal et al.\ \cite{mittal} point out that methods proposed for solving stochastic models often use sophisticated subroutines, such as submodular function minimization and Monte-Carlo techniques, that may not always be practical. Robust optimization (e.g., \cite{bert,bertsim,tal}) offers an alternative framework to address uncertainty that avoids distributional assumptions. Instead, it uses deterministic uncertainty sets and performs a worst-case analysis w.r.t.\ to the uncertainty. This addresses the problems arising out of insufficient data and often these models are more tractable.\footnote{Sometimes at the expense of being too conservative.} Indeed, for the case of appointment scheduling, Mittal et al.\ \cite{mittal} introduced the following robust formulation for the problem.   
	
	Referring to the procedures/services in any given context simply as `jobs', in the robust appointment scheduling problem we are given $n$ jobs with uncertain service times $\{t_i\}$. Assume that the time $t_i$ for job $i$ can be anything in the range $[p_i-\hat{\delta}_i, p_i + \delta_i]$. The range is known to us but the value of $t_i$ can be chosen adversarially. Our task is to propose an appointment start time for each job. A job will be available to process after this start time and jobs will be served in order of increasing index $i$. The order in which jobs are processed is fixed a priori. Due to uncertainty in service times, a job may have to wait to be processed after its appointed start time. The overage cost incurred per unit of wait is given by $o_i \ge 0$. Similarly,  $u_i \ge 0$ is the per-unit underage cost associated with job $i$. Given a set of appointment start times, an adversary chooses the worst possible instantiation of service durations $t_i$ i.e., one that maximizes the overall overage/underage cost. We seek start times to minimize this worst case total cost.  
	
	Mittal et al.\ \cite{mittal} found closed-form optimal solutions for the special case of homogeneous underage costs i.e., $u_i=u$ for all jobs $i$. We study more general cases\footnote{Arguably, just as processing times and uncertainty vary across procedures, so can personnel and equipment. Thus, in general the per-unit underage costs will vary across procedures. Moreover, it is also reasonable to assume that the per-unit cost changes with the amount of underage/overage. For instance, keeping patients waiting for another unit of time becomes less and less desirable with every unit of delay. This prompts the study of piece-wise linear overage and underage costs.}, where a direct generalization of their solution can be arbitrarily bad. Our central goal here is to find simple, easy to implement, and theoretically well founded algorithms. In the following we summarize our contributions and discuss closely related work.

\textbf{Robust Appointment Scheduling (RAS):} 
We explore the structure of optimal appointment schedules and find several properties that all optimal schedules must satisfy. Somewhat surprisingly, we find that in every optimal appointment schedule, the case of all jobs underaged as well as the case of all jobs overaged are always worst-cases for the adversary's problem. Casting these in the form of linear (in)equalities gives an LP, resulting in a 2-approximation. Further simplification yields a closed-form solution to the LP that could be implemented even on a spreadsheet. More generally, when the per-unit costs are allowed to change with the amount of overage/underage and the cost functions are piece-wise linear and convex (and non-decreasing), we are able to generalize these properties to obtain a 2-approximation for the problem via solving a min-cost flow problem on a graph with convex piece-wise linear arc costs. 

Previously, Mittal et al.\ \cite{mittal} examined the special case of identical underage costs and, quite remarkably, found a closed-form optimal solution in that case. While their analysis was tailored for the special case, our different approach allows us to characterize exactly when such a result holds, and we discover more general conditions under which there exists a closed-form optimal solution to the problem. We also show a similar extension for the case of convex piece-wise linear costs, where we give a simple iterative algorithm that gives an optimal allocation under additional assumptions on the cost functions.  

\textbf{Robust Appointment Scheduling and Ordering (RASO):}
While the order of patients is often fixed in advance, there are instances where patient order is flexible and simultaneously part of the decision making \cite{cayi}. When underage costs are identical, Mittal et al.\ \cite{mittal} established a key connection between RASO and the theory of scheduling jobs on a single machine to minimize the sum of a weighted nonlinear (concave) cost function of the completion times \cite{megow,smith,eptas}. Focusing on the special case, we exploit this connection further to give an extremely practical ratio based ordering policy, inspired by Smith's rule \cite{smitho}.  We call it the Customized-Smith rule (short C-Smith) and show that it is 1.0604-approximate. 
Further, we find an algorithm that is as good as C-Smith on every instance, but is also locally optimal w.r.t.\ pairwise swaps of neighboring jobs. Previously, Mittal et al.\ \cite{mittal} showed that using the EPTAS for the related min-sum scheduling problem given by Stiller and Wiese \cite{eptas}, yields an EPTAS for RASO. The runtime of the EPTAS scales as $O(2^{1/\epsilon^7}n)$, and its implementation complexity makes it hard to use in practice. Mittal et al.\ \cite{mittal} offered Smith's rule as a practical alternative, with approximation ratio 1.137 due to a result in H{\"o}hn and Jacobs \cite{smith}. 
Finally, for the general case (arbitrary underage costs), we show the problem is strongly NP hard. We also briefly discuss the inefficacy of \emph{list ordering} heuristics for the problem and give a heuristic with a matching lower bound.

\textbf{Other related work:}
The existing literature on appointment scheduling is quite diverse and, in addition to the robust model we discussed above, includes numerous results on stochastic optimization models, queueing models, as well as distributionally robust optimization models. We only give a very brief review here and refer the reader to \cite{cayi} for a broad overview of the challenges in scheduling appointment systems in health care, and to \cite{mak,mittal} for more comprehensive surveys on past work. 
Starting with more recent work, Jiang et al.\ \cite{jiang} consider a distributionally robust model and propose a MINLP formulation that can handle random no-shows. Mak et al.\ \cite{mak} consider a distributionally robust model 
given marginal moments for job durations, and formulate the problem as a tractable conic program. For the case of flexible job order, they show that under certain assumptions a simple variance based ordering can be optimal. 
Prior to \cite{mak}, Kong et al.\ \cite{kong} first considered a distributionally robust formulation but with cross-moments as opposed to marginal moments. They formulate the problem as a copostive cone program and propose a tractable SDP relaxation. Wang \cite{wang1} modeled the problem using a queueing model where the processing times
of the jobs are i.i.d.\ exponential and new jobs may be released over time. Wang \cite{wang2} generalized the model to allow for different mean processing times for jobs. 

In a different direction, the problem was modeled as a two-stage stochastic linear program in Denton and Gupta \cite{denton}. 
For this problem, Robinson and Chen \cite{robin} compute near-optimal solutions using a Monte-Carlo integration technique.  
Kaandorp and Koole \cite{kaand} considered a local search algorithm and proved that it converges to an optimal solution. Another stochastic model was introduced by Green et al.\  \cite{green}. They considered the problem of outpatient appointment scheduling with emergency services and modeled it as a dynamic stochastic control problem. In highly influential work, Begen and Queyranne \cite{begen} considered a discrete stochastic model (where job durations are integer random variables with finite support). They showed that the problem reduces to an instance of submodular function minimization, under certain assumptions on the per-unit costs. Begen et al.\ \cite{begens} extended the idea, proving a similar result for a data-driven discrete stochastic model. Ge et al.\ \cite{genote} further extended the result of Begen and Queyranne \cite{begen}, to the setting of convex piece-wise linear per-unit costs. 

Next, we briefly discuss some of the assumptions on cost coefficients made in previous work, and how our assumptions fit into the existing picture.

 \textbf{Assumptions on underage costs:} Previous work on the stochastic and distributionally robust models discussed above, often assumes an upper bound on the variation in underage costs. Formally, consider jobs that are indexed in the order they are scheduled. So we have jobs $i\in\{1,\dots,n\}$,with per unit costs $u_i,o_i$. Now for instance, the result by Begen and Queyranne \cite{begen}, assumes that there exists values $\alpha_i\leq o_i$ such that $u_{i+1}\leq u_{i}+\alpha_i -\alpha_{i+1}$ for every $i\in[n-1]$. Descending underage costs for example, satisfy this condition. More generally, \cite{jiang,mak,kong,genote} assume that $u_{i+1}\leq o_i+u_i$ for every $i\in[n-1]$ (or that $u_i=u$ for all $i$). In effect, these assumptions are satisfied as long as the underage costs don't \emph{increase} too sharply.
 
 In contrast, we give a closed form optimal allocation for instances where the $u_i$ values do not \emph{decrease} too sharply. For example, our results yield a closed form optimal solution for the case of non-decreasing underage costs. This includes the case $u_{i+1}\geq u_i+o_i\, \forall i\in[n-1]$, converse of the condition assumed in some previous work. If the per-unit underage costs are decreasing, we show that as long as the per-unit values are ``large enough'', our LP still leads to a closed form optimal allocation. In practical terms, it is quite possible that the underage costs could be increasing or decreasing, preventing a direct comparison of our assumptions with past work. In fact, for the RASO problem where one has the flexibility to choose the order of jobs, an optimal order can have arbitrarily varying underage costs. From a technical viewpoint, our results complement previous work, indicating that perhaps the robust model is more tractable in instances where other models are not, and vice versa. 
 
 \textbf{Overview of the paper:} In Section \ref{notation}, we introduce notation and formally state the problem. We also make certain observations that simplify the problem w.l.o.g..\ Then in Section \ref{RAS}, we explore the appointment scheduling problem under fixed order. We first give our LP based approximation for general underage costs in Section \ref{RASG}. And later in Section \ref{sec:optspec}, we tighten our LP formulation to obtain optimal closed form solutions under additional assumptions on the costs. In Section \ref{sec:rass}, we consider the problem of jointly finding an optimal order as well as appointment schedule (RASO). We show that the problem is NP hard in general and discuss the limitations of a classes of simple heuristics in Section \ref{rasog}. In Section \ref{rasohomo}, we focus on the ordering and scheduling problem in the special case of homogeneous underage costs. Here we discuss some disadvantages of existing results and propose two new heuristics to tackle some of the issues. 
 Finally, we conclude with some open problems in Section \ref{conclusion}. 


\section{Notation \& Preliminaries}\label{notation}
We start with a description of RAS, and define additional notation required for RASO in Section \ref{sec:rass}. Recall, we have $n$ jobs to be served in order of increasing index $i\in \{1,\dots,n\}$. Service times are uncertain and modeled via a box uncertainty set: job $i$ takes time $t_i$ in the range $[p_i-\hat{\delta}_i, p_i + \delta_i]$.  Here we assume that $0 \le \hat{\delta}_i \le p_i$ as well as $\delta_i \ge 0$, for all jobs $i$.  We denote by $o_i \ge 0$ and $u_i \ge 0$ the per-unit overage cost and per-unit underage cost of job $i$, respectively. More generally, we represent overage and underage costs as \emph{non-decreasing} functions $o_i(\cdot),u_i(\cdot)$ respectively. The case of constant per-unit costs is then given by $o_i(x)=o_i\cdot x$ and $u_i(x)=u_i\cdot x$. We would like to appoint job start times, $\{A_i\}_i$, such that every job arrives at its appointed time and is served as soon as possible after.  
Given $\{A_i\}_i$, consider an arbitrary instance of service times $\{t_i\}_i$, and let $\{C_i\}_i$ denote the completion times of jobs.\footnote{ $C_i=\max\{C_{i-1},A_i\}+t_i$ for all $i\geq 1$ with $C_0=A_1=0$, w.l.o.g.} If job $i$ is delayed and ends after the appointed start time for job $i+1$, we incur overage cost $o_i (C_i-A_{i+1})$. Similarly, if job $i$ ends before the appointed start time for job $i+1$, we incur underage cost $u_i(A_{i+1}-C_i)$. Therefore, the cost due to job $i$ is $\max\{o_i(C_i-A_{i+1}),u_i(A_{i+1}-C_i)\}$. 
The RAS problem can now be stated as follows,
\begin{equation}
\min_{\{A_i\}_i}\quad \underbrace{\Big(\max_{t_i\in [p_i-\hat{\delta}_i, p_i + \delta_i]\,  \forall i\in [n]}\quad \sum_{i=1}^n \max\{o_i(C_i-A_{i+1}),u_i(A_{i+1}-C_i)\}\Big)}_\textit{Adversary's problem, given appointment times $A_i$.}. \label{problem}
\end{equation}
The adversary's problem in (\ref{problem}), finds a \textit{worst possible profile/instance of service times} $T=(t_i)=(t_1,\dots,t_n)$, maximizing the cost, given the schedule $\{A_i\}_i$. Let $c(T,A,i)= \max\{o_i(C_i(T)-A_{i+1}),u_i(A_{i+1}-C_i(T))\}$ denote the cost of job~$i$ given allocation $\{A_i\}_i$ and time profile $T$. When the allocation is clear from the context, we use the shorthand $c(T,i)$. Let $c(T)=\sum_i c(T,i)$ denote the total cost due to profile $T$. It is not difficult to see that an optimal allocation needs to allocate at least $p_i-\hat{\delta_i}$ time for job $i$, so we let $p_i-\hat{\delta_i}=0$ for every job $i$, w.l.o.g.\ (also in Lemma 6 of \cite{mittal}). To simplify notation henceforth, we let service times be in the range $[0,\Delta_i]$, where $\Delta_i=\hat{\delta_i}+\delta_i$. Now, given appointment schedule $\{A_i\}_i$, consider equivalent variables $x_i$ that represent the duration allocated for a job. We have, $x_i=A_{i+1}-A_{i}$, which is the difference between the start times of job $i+1$ and job $i$. Equivalently, $A_i=\sum_{j=1}^{i-1} x_j$ for $i\geq 2$ and $A_1=0$, w.l.o.g. We call $\{x_i\}_i$, and sometimes by abuse of terminology $\{A_i\}_i$, the \textit{allocation}. 
Note that in our model, job $n$ is cost free since there is no appointment succeeding it ($A_{n+1}=\infty$). W.l.o.g.\ we may assume that job $n$ also suffers from overage and underage based on its assigned end time (assign a dummy job to succeed job $n$). Note, if $u_i(\cdot)=0$, we allot a very large time slot for job $i$ and jobs on different sides of the slot become independent. Therefore, we let\ $u_i(\cdot)\neq 0$ for all $i$, w.l.o.g..\ Also, if $o_n(\cdot)=0$ we can assume that $x_n=0$ and in fact ignore job $n$, therefore we also let $o_n(\cdot)\neq 0$, w.l.o.g..\ To coalesce, the assumptions we have made w.l.o.g., so far: 

\emph{(i) Service time for job $i$ lies in $[0,\Delta_i]$; (ii)  for every $i$, $u_i(\cdot)\neq 0$  ($u_i>0$ in the constant per-unit case); and (iii) $o_n(\cdot)\neq 0$ ($o_n>0$ in the constant per-unit case).}

Finally, we say job $i$ is \emph{underaged} if it ends on or before time $A_{i+1}$, and \emph{overaged} otherwise. However, if $\Delta_i=0$ and job $i$ starts/ends at $A_{i+1}$, we consider it to be both overaged and underaged (for technical reasons). 
Further, if job $i$ takes time $\Delta_i$ to be served, we say it runs for \emph{maximum time}. Note that if a job $i$ is underaged in some worst-case $T=(t_i)$, then w.l.o.g., $t_i=0$, as the underaged costs $u_i(\cdot)$ are non-decreasing functions for every $i$ (similar to Lemma 11 of \cite{mittal}). Observe also that if a job $i$ is underaged in some worst-case time profile, and it has zero cost, then $i$ must end at $A_{i+1}$ and take zero time. Otherwise, we have a strictly worse case by underaging $i$ (simply reduce the time taken by $i$ by some small $\epsilon>0$).  

Finally, let $S_i=\{i,i+1,\dots,n\}=\{j| j\geq i\}$ be the subset of the last $n-i+1$ jobs in the schedule. Let $\x{S}{i}$ denote the optimal time allocation for job $i$ when considering only jobs in a subset $S$ that contains $i$. When $S=[n]$ we often use shorthand $x_i$. The next section discusses a result from previous work and offers some intuition for our treatment of the general case that follows subsequently in Section \ref{RAS}.

\subsection{Closed Form Optimal Solution of Mittal et al. \cite{mittal}}

Consider a single job, with per unit costs $u,o$ and maximum time $\Delta$. If this job is alloted time duration $x$, then the worst case cost is given by $\max\{u(x),o(\Delta-x)\}$, minimized at $x=\Delta \frac{o}{o+u}$ (recall, we assume that a dummy job always follows the last job, so the last job incurs overage cost for delays). Mittal et al.\ \cite{mittal} showed that for constant per-unit underage costs, $u_i=u$ for every job $i\in[n]$, this formula generalizes and the optimal allocation is given simply by, $x_i=\Delta_i\frac{ o(i)}{o(i)+u}$. Here $o(i)$ is the sum $\sum_{j=i}^{n}o_j$, of the per unit overage costs of the jobs succeeding job $i$. So each job is alloted a fraction of its maximum service time $\Delta_i$, with a smaller fraction for larger values of $u$ (to prevent large underage costs). Further, earlier jobs are alloted a larger fraction of their maximum time to prevent a large cascade of delays for jobs serviced later on. 

Now, consider a natural generalization of this formula for heterogeneous underage costs given by, $x_i=\Delta_i\frac{ o(i)}{o(i)+u_i}$. Unfortunately, this can be a suboptimal allocation even for two jobs with reasonable parameter values, and in general an arbitrarily bad approximation. For example, consider two jobs that are almost identical except that job 2 (which is scheduled later) has a small underage cost. Specifically, let $\Delta_1=\Delta_2=1$ and suppose per unit costs $o_1=o_2=1$ and $u_1=M\geq 2$, but $u_2=1$. The allocation given by the formula sets $x_1$ to $\frac{2}{M+2}$ and $x_2 $ to $\frac{1}{2}$. The worst case cost of this allocation is attained when both jobs are underaged or both overaged (more on this later), and equals $\frac{2M}{M+2}+\frac{1}{2}$. For large $M$, this value approaches $\frac{5}{2}$. Now instead, consider the following allocation, $y_1=0$ and $y_2=\Delta_2\frac{o_2}{o_2+u_2}+\Delta_1\frac{o(1)}{o_2+u_2}=\frac{1}{2}+1$ (we show later in Section \ref{sec:optspec}, that this is in fact an optimal allocation). It is easily checked that the worst case cost is $\frac{3}{2}$ for this allocation. More generally, if we set $o_1$ and $u_2$ to a small value $\epsilon$, the first allocation becomes arbitrarily bad. Intuitively, this demonstrates that if jobs that are later in the order have small underage costs, it is beneficial to allocate larger time to these later jobs to buffer for delays from earlier jobs. When this is not the case (such as for homogeneous underage costs), it is better to instead buffer by allocating larger fractions to jobs that are earlier in the order.  


\section{Robust Appointment Scheduling}\label{RAS}
\subsection{Heterogeneous Per-unit Underage Costs}\label{RASG}
In this section, we focus on the case of constant per-unit costs $o_i,u_i$. The key behind our results lies in finding useful properties satisfied by optimal allocations. Towards that end, it will be instrumental to understand the worst cases (solutions to the adversaries' problem) for optimal allocations. Recall the adversary's objective: given $\{A_i\}_i$, maximize cost $c(T,A)$ over all possible time profiles $T$ allowed by the uncertainty set. 
  It turns out that for an arbitrary allocation, there could be a unique worst-case where some jobs are overaged while others are underaged, and for some jobs the service time $t_i$ is neither $0$ nor $\Delta_i$ \footnote{This is in contrast to Lemma 9 in \cite{mittal} for the special case of $u_i=u\, \forall i$.}. It turns out though, that for any given allocation the adversary's problem can be reduced to an instance of finding longest paths on a directed acyclic graph with $n+1$ nodes \cite{jim}, and therefore can be solved in polynomial time. However, using this fact to find an optimal allocation does not obviously lead to a tractable problem. 
Instead, here  we shall find properties in the form of linear inequalities that every optimal allocation satisfies. Using these, we formulate an LP relaxation for the problem and show that an optimal solution to the LP is 2-approximate. 

As a natural next step, we then look for more structure to further tighten the formulation. Surprisingly, we find that the adversary's problem given an \emph{optimal allocation} for problem (\ref{problem}) is easily solved; in every optimal allocation, $t_i=\Delta_i$ for every $i$ (all jobs taking maximum time) and $t_i=0$ for all $i$ (all jobs taking zero time) are worst-cases. Combining this with other structural insights, in Section \ref{sec:optspec} we propose a strengthened formulation that leads to closed form optimal solutions under some assumptions on underage costs. In the lemmas that follow we introduce two properties of optimal allocations, leading to the first LP formulation.  


\begin{lem}\label{sanity1}
	For every optimal allocation $\{x_i\}_i$, $\sum_{j=1}^{k}x_{j}\leq \sum_{j=1}^{k} \Delta_j$\ for all $k = 1, \ldots, n$. 
\end{lem}
\begin{proof}
	Given an optimal allocation $\{x_i\}_i$, let $\delta_k=\sum_{j=1}^{k}x_{j}- \sum_{j=1}^{k} \Delta_j$, for $k\in[n]$. Suppose, $\delta_k>0$ for some $k$, and let $k_0$ be the smallest such $k$. Given this, notice that even if the first $k_0$ jobs take maximum time, job $k_0$ can never be overaged. Therefore, decreasing $x_{k_0}$ decreases the underage cost of job $k_0$ ($u_{k_0}>0$) and thus, also the worst-case cost. So the allocation $\{x'_i\}_i$ where, $x'_{k_0}=x_{k_0}-\delta_{k_0}$ and $x'_i=x_i$ for all other $i$ is clearly a better allocation, contradiction. 
	This lemma also follows as a direct corollary of Lemma \ref{l1} stated later on. 
\end{proof}
\begin{lem}\label{costub}
	Given an allocation $\{x_i\}_{i}$ where $\sum_{j=1}^{i}x_{j}\leq \sum_{j=1}^{i} \Delta_j$\ for all $i$, and a time profile $T=(t_i)$. The cost $ c(T)$ of profile $T$ is at most $\sum_{i=1}^n o(i)(\Delta_i-x_i)+\sum_{i=1}^n u_i x_i$.  
\end{lem}
\begin{proof}
	For any given time profile $T$, job $i$ is either underaged with $c(T,i)\leq u_ix_i$ or overaged with cost at most $o_i\sum_{j \leq i} (\Delta_j-x_j)$. The latter follows from the fact that $C_i\leq \sum_{j \leq i}\Delta_j$ since $\sum_{j=1}^{k}x_{j}\leq \sum_{j=1}^{k} \Delta_j$ for all $k\in[i]$. Therefore, $\sum_{i=1}^n c(T,i)\leq \sum_{i=1}^n o_i\sum_{j=1}^{i} (\Delta_j - x_j)+\sum_{i=1}^n u_i x_i$. Rearranging the sum we have, $\sum_{i=1}^n o_i\sum_{j=1}^{i} (\Delta_j - x_j)= \sum_{j=1}^n (\Delta_j - x_j)o(j)$.
\end{proof}

Using the above results, the following LP now gives us an extremely simple approximation for the general problem that could be easily implemented on most systems. 
\begin{eqnarray}
\textbf{LP-1:}\quad  &\min &\sum_{j=1}^{n} (u_j-o(j))y_j \nonumber\\
&{\rm s.t.} 
 & \sum_{j=1}^{k}(y_j-\Delta_j)\leq 0\quad \forall k\in [n];\quad \label{aux}\\
&& y_j\geq 0 \quad \forall j\in [n].\nonumber
\end{eqnarray} 

\begin{thm}\label{closed2}\label{2apx}
	LP-1 is a (tight) $2$-approximation, and the following is an optimal solution to LP-1. Define $m_i=\min\big(\arg\min_{j\geq i}(u_j-o(j))\big)$, then for every $i$, 
	\[
  y_i =
  \begin{cases}
    0 & \text{if $u_i-o(i)\geq 0$ or $i\neq m_i$} \\
    \sum_{j|m_j=i}\Delta_j & \text{otherwise.}
  \end{cases}
\]
\end{thm}
\begin{proof}
	The proof of optimality for the proposed solution is easy to verify. To see the guarantee let us re-write the LP-1 objective as, \[\min \sum_ju_jy_j +\sum_j o(j)(\Delta_j-y_j).\] From constraints~\eqref{aux} and Lemma \ref{costub}, this is an upper bound on the worst-case cost of any feasible solution to LP-1. Let $\{y^*_i\}$ denote an optimal solution to the LP. Lemma \ref{sanity1} implies that every optimal allocation $\{x_i\}$ is a feasible solution to LP-1. Further, $\sum_ju_jx_j$ denotes the cost when all jobs take zero time and $\sum_j o(j)(\Delta_j-x_j)$ is the cost of overaging all jobs with each job taking maximum time ($\Delta_i$). Therefore, the worst-case cost of an optimal allocation, denoted $OPT$, is at least $\max\{\sum_ju_jx_j,\sum_j o(j)(\Delta_j-x_j)\}$ (we show later that these two costs are both in fact, equal to $OPT$). 
	Therefore, \[ \sum_ju_jy^*_j +\sum_j o(j)(\Delta_j-y^*_j)\leq \sum_ju_jx_j +\sum_j o(j)(\Delta_j-x_j)\leq 2~OPT.\]
	
	For a tight instance, consider two jobs $\{1,2\}$ with $\Delta_1=\epsilon\to 0$, $u_1=1/\epsilon, o_1=1/\epsilon-(1+\epsilon)$ and $\Delta_2=u_2=o_2=1$. Therefore, $\frac{u_2}{u_2+o_2}=0.5\lessapprox \frac{u_1}{u_1+o(1)}$. Consider the allocation $x_2=0.5$ and $x_1=\epsilon\frac{o(1)}{u_1+o(1)}\approx \epsilon/2$. It is easy to see that the worst-case of the allocation is when both jobs are underaged (overaged) and hence the cost of this allocation is $u_1x_1+u_2x_2\approx 1$. Now, consider the solution $y_2=0.5+\frac{\Delta_1o(1)}{u_2+o_2}\approx 1$ and $y_1=0$. This is an optimal solution to the LP and the worst-case occurs when job 2 is underaged and job 1 is overaged. The worst-case cost of this allocation is $ o_1\Delta_1+ u_2(y_2-\Delta_1)\approx 2$. Notice that in this example, the formula allocation $x_i$ is near optimal while the optimal solution to the LP is not (here $u_2-o(2)=0$, so LP-1 has multiple optimal solutions, and the solution $y_1=y_2=0$ given by the closed form above also has worst-case cost $\approx 2$). However, while the formula allocation can be arbitrarily bad in general, the LP optimal is guaranteed to be within a factor of 2.
\end{proof}
\subsubsection{Beyond Constant Per-unit Costs} \label{2apxg}
Suppose instead of scalar costs $o_i,u_i$, we have non-decreasing, piece-wise linear and convex cost functions $o_i(\cdot),u_i(\cdot)$. Then it is easy to check that there exists an optimal allocation satisfying the conditions in Lemma \ref{sanity1} (see Lemma \ref{sanity1g} in Appendix \ref{ext1}). Similar to Lemma \ref{costub}, we have that given an allocation $\{x_i\}_i$, satisfying Lemma \ref{sanity1}, the cost $c(T)$ for any profile $T$ is at most, \[\sum_i u_i(x_i) + \sum_i o_i\Big(\sum_{j\leq i} (\Delta_j-x_j)\Big).\] 
 Now the problem of minimizing this objective subject to the linear constraints in LP-1, is a min cost flow problem with arc costs given by the overage and underage cost functions. More specifically, consider a directed graph with $n+1$ nodes and $2n$ edges, where there is a directed edge from node $i$ to node $n+1$ with cost $u_i(\cdot)$ and an edge $i \to i+1$ with cost $o_i(\cdot)$,  for every $i\in[n-1]$. Finally, there is another edge from $n$ to $n+1$ with cost $o_n(\cdot)$. Finally, each node $i$ has a supply of $\Delta_i$ and node $n+1$ is a sink with demand $\sum_i\Delta_i$. Given a feasible flow in this graph, the flow on edge $(i\to n+1)$ gives the time allocation for job $i$, and vice versa. Now, if the costs are piece-wise linear and convex, we have from the algorithm by Pinto and Shamir \cite{pinto} for min flows with convex piece-wise linear costs, that the problem can be solved in polynomial (in $n$ and the maximum number of pieces in the cost functions) time. The optimal solution to this problem is a 2-approximation, and the analysis closely resembles the proof of Theorem \ref{closed2} (see Appendix \ref{ext1}).
\subsection{Optimal Solution to RAS for Special Cases}\label{sec:optspec}
Let us now investigate additional properties with the goal of strengthening our LP. We start by proving our claim from earlier --  all jobs underaged (taking 0 time) and all jobs overaged (taking maximum time) are worst-cases for optimal allocations. We break down the proof into smaller parts. The first lemma is very useful and appears often in proofs of other lemmas. The key insight behind the lemma is simple -- when all jobs overaged is a worst-case, if the first job is forced to start late, then the case of all jobs overaged suffers maximum increase in cost, and is thus still a worst case.

\begin{lem}\label{shift}
	Given an arbitrary allocation $\{A_i\}_i$, where $t_i=\Delta_i$ for all $i$ is a worst-case and every job $i$ completes at or after time $A_{i+1}$ in this case (this holds, for instance, if $\sum_{j\leq i}x_j\leq \sum_{j\leq i}\Delta_j$ for all $i$). Recall that w.l.o.g., $A_1=0$ and consider a modified problem for the adversary, where the first job is always forced to start at a later time $t_0$ instead of time $A_1=0$. Then, $T=(t_i)$ is also a worst-case for the modified problem.
\end{lem}
\begin{proof}
	Let $T=(\Delta_i)$ denote the profile for all jobs taking maximum time. When job $1$ starts at time 0, denote the cost of job $i$ by $c(0,T,i)$. Since $T$ is a worst-case by assumption, we have $c(0,T)\geq c(0,Z)$ for every profile $Z$. For the modified setting where job $1$ starts at time $t>0$, we claim that $c(t,Z,i)\leq c(0,Z,i)+o_it$ for every $Z$. To see this, suppose job 1 starts at 0 and consider two cases: (i) $i$ overaged in $Z$ and (ii) $i$ underaged in $Z$. In case (i), $i$ will still be overaged when job 1 starts at time $t$ and the completion time of $i$ can increase by at most $t$. In case (ii), if $i$ is still underaged when job 1 starts at time $t$, we are done. Else, $i$ becomes overaged but the maximum overage cost is $o_it$. Now for profile $T$, since every job $i$ completes at or after $A_{i+1}$ when job 1 starts at time 0, $c(t,T,i)=c(0,T,i)+o_it$. Therefore, $c(t,Z)\leq \sum_i (c(0,Z,i)+o_i t)\leq \sum_i c(0,T,i)+\sum_i o_i t= c(t,T)$.
\end{proof}

The next lemma says that we cannot have an optimal allocation where a certain job is underaged (overaged) in all the worst-case time profiles. Observe that if  job $i$ is always underaged, simply reducing the allocation $x_i$ would give a strictly better allocation, contradicting optimality. Indeed, Lemmas 12 and 15 in \cite{mittal} argue exactly this. However, that argument fails if there exists a worst-case where $i$ is underaged with zero cost (occurs when $i-1$ ends at $A_{i+1}$ and $t_i=0$, which we defined as a case of underage in Section \ref{notation}). This gives rise to a subtle issue that demands a more involved argument (similarly for the case of overage). We postpone the formal proof to Appendix \ref{apxoptspec}. 
\begin{lem}\label{l1}
	Suppose $o_j>0$ for every job $j$. Given an optimal allocation, consider an arbitrary job $i$. There exists a worst-case $T=(t_i)$ where $i$ is underaged and $t_i=0$, as well as a worst-case where $i$ is overaged (with some $t_i$ that is not necessarily $\Delta_i$). 
\end{lem}

The assumption $o_j>0$ is only needed in the very last case of the proof. If $o_j=0$ for some jobs $j$, a simple limiting argument (also in Appendix \ref{apxoptspec}) shows that Lemmas \ref{worst2} and \ref{worst1} below still hold for at least one optimal allocation, and this is all that we use in subsequent results.
\begin{lem}\label{worst2}
	In every optimal allocation (in at least one optimal allocation, if $o_j=0$ for some $j$), the case of all jobs underaged, i.e.,\ $t_i=0$ for every job $i$, is a worst-case.
\end{lem}
\begin{proof}
	Let $\{A_i\}_i$ denote an optimal allocation. Lemma \ref{l1} implies there is a worst-case where job $1$ takes zero time. We proceed via induction, assuming there is a worst-case $T=(t_i)$ where $t_i=0$ for $i\in[k]=\{1,\dots,k\}$, i.e., all jobs in $[k]$ are underaged. We will show there exists a worst-case where all jobs in $[k+1]$ are underaged. 
	
	 Suppose that $k+1$ is overaged in $T$ (otherwise we are done). Lemma \ref{l1} implies there is a worst-case, denoted $T'=(t'_i)$, where job $k+1$ is underaged and takes zero time. Let $C_{k+1}$ denote the completion time of job $k+1$ in $T$. Clearly, $C_{k+1}> A_{k+2}$ and $k+1$ starts at $A_{k+1}$ in $T$. Similarly, let $C'_{k+1}$ denote the start/completion time of job $k+1$ in $T'$. Then, $C'_{k+1}\leq A_{k+2}$. Now consider a new profile $Z=(z_i)$, formed by a combination of $T$ and $T'$. We let $z_i=t'_i$ for $i\in[k]$ and $z_i=t_i$ for $i\in\{k+2,\dots,n\}$. Since the completion times of job $k$ are identical in $T'$ and $Z$, we set $z_{k+1}=C_{k+1}-C'_{k+1}\leq C_{k+1}-A_{k+1}=t_{k+1}$. Therefore, job $k+1$ ends at time $C_{k+1}$ in $Z$. Now, 
	observe that $\sum_{i\in [n]} c(Z,i)=\sum_{i=1}^{k} c(T',i) + \sum_{i=k+1}^{n}c(T,i)$. Since $T$ is a worst-case, we also have $\sum_{i=1}^{k} c(T',i) + \sum_{i=k+1}^{n}c(T,i)\leq \sum_{i=1}^{n}c(T,i)$. Therefore, 
	$\sum_{i=1}^{k} c(T',i)\leq \sum_{i=1}^{k}c(T,i).$
	Now, consider another hybrid case $Q$, where jobs $1$ to $k$ are all underaged and take zero time as in $T$ and jobs $k+2$ to $n$ are as in $T'$. Job $k+1$ starts/ends at $A_{k+1}$ in $Q$ and hence $c(T',k+1)\leq c(Q,k+1)$. Combining everything, $\sum_i c(T',i)\leq\sum_{i=1}^{k} c(T,i) + \sum _{k+1}^{n} c(T',i)\leq \sum_{i=1}^{k} c(T,i) + c(Q,k+1)+ \sum _{k+2}^{n} c(T',i) = \sum_i c(Q,i).$
	Hence, $Q$ is a worst-case with jobs $1$ to $k+1$ all underaged.
\end{proof}

\begin{lem}\label{worst1}
	In every optimal allocation (in at least one optimal allocation, if $o_j=0$ for some $j$), the case of all jobs overaged is a worst-case. Therefore, $t_i=\Delta_i$ for all $i$ is a worst-case.
\end{lem}
\begin{proof}
	Let $\{A_i\}_i$ denote an optimal allocation. Observe that if there exists a worst-case where jobs $S_k=\{k,\dots,n\}$ are all overaged for some $k$, then there exists a worst-case where all jobs in $S_k$ are overaged and take maximum time. 
	We proceed by induction and show that if there is a worst-case where every job in $S_{k}$ is overaged, then there exists a worst-case where every job in $S_{k-1}$ is overaged. For $k=n$, by Lemma \ref{l1} there is a worst-case where job $n$ is overaged. Assume there exists a worst-case where jobs $k$ to $n$ are overaged, denoted as $T=(t_i)$. 
	
	 Suppose that job $k-1$ is underaged in profile $T$ (otherwise we are done). Recall that $u_{k-1}>0$ and $t_{k-1}=0$, since $k-1$ is underaged in $T$. Then, since $T$ is a worst-case profile; restricted to the subset $S_{k}$, the profile $\{t_i=\Delta_i\}_{i\geq k}$ is a worst-case with all jobs overaged, for allocation $\{A_i\}_{i\geq k}$. Now, by Lemma \ref{l1} for the set of all jobs $[n]$, there exists a worst-case where job $k-1$ is overaged, denoted $T'=(t'_i)$. Then, consider profile $Z=(z_i)$ with $z_i=t'_i$ for jobs in $S_1-S_{k-1}$ and $z_i=\Delta_i$ for $i\in S_{k-1}$. $Z$ is a worst-case profile since it matches the cost in $T'$ for jobs in $S_1-S_{k-1}$, and due to Lemma \ref{shift} the total costs for jobs in $S_{k-1}$ can only be larger than the same total in $T'$.   
\end{proof}

Therefore, the following equation characterizes the worst-case cost of every optimal allocation $\{x_i\}_i$,
\begin{equation}\label{sumcost}
	 \sum_{i=1}^n o(i)(\Delta_i-x_i)=\sum_{i=1}^n u_i x_i. 
\end{equation}
We can already modify LP-1 by adding the equality above. However, this property alone is not sufficient to offer improved results even for identical underage costs. To that end, we introduce the following technical property.
 
\begin{lem}\label{l3}
	Given an optimal allocation $\{x_i\}_{i\in [n]}$ for which all jobs underaged and all jobs overaged are worst-cases (by Lemmas \ref{worst2} and \ref{worst1}, every optimal allocation if $o_j>0$ for all $j$), 
	\begin{equation}
	\sum_{j=k}^{n} (u_j+o(j))x_j\geq \sum_{j=k}^{n} o(j)\Delta_j\quad  \text{for all } k\in[n]. \label{knap}
	\end{equation}
\end{lem}
To get some intuition behind the lemma, consider the objective of minimizing $\sum_i u_iy_i$ for non-negative $y_i$ subject to the equation~\eqref{sumcost}. A greedy solution that sets $y_i=\frac{\sum_j \Delta_jo(j)}{u_i+o(i)}$ for $i=\arg\min_{j\in[n]}\frac{u_j}{u_j+o(j)}$ is optimal. Here the job with the minimum ratio $\frac{u_i}{u_i+o(i)}$, bears the entire burden of the equality~\eqref{sumcost}. Lemma \ref{l3} says that this can only occur when $i=n$ and more generally, places a lower bound on the contribution to equality~\eqref{sumcost} from values $y_k$ to $y_n$, for all $k$. The proof of the lemma is rather technical and is postponed to Appendix \ref{apxoptspec}, along with two accompanying helper lemmas. 
Consider now the strengthened formulation,
\begin{eqnarray}
\textbf{LP-2:}\quad &\min &\sum_{j=1}^{n} u_jy_j\nonumber\\
&{\rm s.t.} &\sum_{j=1}^{n} (u_j+o(j))y_j= \sum_{j=1}^{n} o(j)\Delta_j\label{bala}\\
&& \sum_{j=k}^{n} (u_j+o(j))y_j\geq \sum_{j=k}^{n} o(j)\Delta_j\quad  \forall k\in[n] \label{techineq}\\
&& \sum_{j=1}^{k}(y_j-\Delta_j)\leq 0\quad \forall k\in [n]\nonumber\\
&& y_j\geq 0 \quad \forall j\in [n]\nonumber
\end{eqnarray} 
Clearly, every optimal allocation (at least one optimal allocation, if $o_j=0$ for some $j$) is a feasible solution for the LP. However, an optimal solution to the LP need not have all jobs underaged (or overaged) as worst-case, and hence need not be an optimal allocation. This is demonstrated by the example used in Theorem \ref{2apx}, where it is easily checked that the additional constraints given by \eqref{bala} and \eqref{techineq} do not improve the worst-case approximation bound from earlier. 
 However, we show that under additional assumptions \text{LP-2} yields optimal allocations.
 
 The recipe behind proving this is as follows: (i)~Using the assumptions, show that there is actually a closed-form optimal solution to the LP. (ii) Show that for this solution, all jobs underaged and all jobs overaged are worst-cases. In particular, using this recipe for the special case of homogeneous underage costs $u_i=u$ for every $i$, we find that $x_i=\frac{o(i)\Delta_i}{u+o(i)}$ for all $i\in[n]$ is an optimal solution to LP-2 and an optimal allocation. More generally, we have the following. 
\begin{thm}\label{genmit}
	If $u_{i}\leq u_{i+1}\frac{o(i)}{o(i+1)}$ for all $i\in[n-1]$, then the allocation given by $x_i=\frac{o(i)\Delta_i}{u_i+o(i)}$, $i\in[n]$, is an optimal allocation. 
\end{thm}	
 This generalizes and offers a different perspective on Theorem 5 in \cite{mittal} (since $u_i=u$ for every $i$, implies $u_{i}\leq u_{i+1}\frac{o(i)}{o(i+1)}$ for all $i\in[n-1]$). As a direct implication of the above, we have a closed form optimal allocation if the underage costs are non-decreasing i.e.,
 \[u_i\leq u_{i+1}\quad \forall i\in[n-1].\] 
 This includes for instance, the case of increasing underage costs where $u_{i+1}\geq u_i+o_i$. 
 As we discussed earlier in Section \ref{intro}, this complements the assumptions made in previous work, where the common assumption is that the underage costs are not increasing too drastically and in particular that, $u_{i+1}\leq u_i+o_i$ for all $i\in[n-1]$.
 
  Next, we claim (proof in Appendix \ref{apxoptspec}) that even if $u_{i}> u_{i+1}\frac{o(i)}{o(i+1)}$ for some $i$, there may still exist an LP optimal solution that is also an optimal allocation. This generalizes Theorem \ref{genmit}. A direct corollary of the result is that if underage costs are ``large enough'', we have a closed form optimal solution even if the costs otherwise vary arbitrarily. More specifically, if the underage costs are such that for every pair of jobs $i-l$ and $i$ with $u_{i-l}>u_i$, we have,
  \[u_{i}\geq \sum_{j=1}^l o_{i-j}.\]
  Then LP-2 leads to a closed form optimal allocation.
  
 \begin{thm}\label{sidedspec}
 	Given $n$ jobs with parameters such that for $i\in[n]$, whenever $i\neq m_i\coloneqq \min\big( \arg \min_{j\geq i} \frac{u_j}{u_j+o(j)}\big)$, we have $u_{m_i}\geq o(i)-o(m_i)$, then the following is both an optimal solution to the LP and an optimal allocation,
 	\[ x_k=\frac{\sum_{i|k=m_i} \Delta_i o(i)}{u_k+o(k)} \quad \text{for all } k\in[n]. \]
 \end{thm}
\subsubsection{Optimal Allocation Beyond Constant Per-unit Costs}\label{optg}
For non-decreasing, piece-wise linear and convex costs $o_i(\cdot),u_i(\cdot)$, let $\olo_i ,\ulo_i$ denote the largest and smallest slope for $o_i(\cdot)$  and, $\olu_i,\ulu_i$ the largest and smallest slope for $u_i(\cdot)$. By generalizing insights from the case of constant per-unit costs, in Appendix \ref{ext2} we show that Algorithm \ref{optgalg} finds the optimal allocation when  $\ulu_{i+1}\geq \olu_i, \text{ for all $i\in [n-1]$.}$ This condition is satisfied for instance if $u_i(x)=u\cdot x$ for all $i\in[n]$, and $o_i(\cdot)$ is an arbitrary non decreasing, convex piece-wise linear function. Similar to the case of constant per-unit costs, our assumption complements assumptions made in previous work on the stochastic setting of the problem. For instance, in the stochastic setting Ge et al. \cite{genote}, gave a polynomial time algorithm for non decreasing, piece-wise linear and convex costs when $\olu_{i+1}\leq \ulu_i+\ulo_i, $ for all $i\in [n-1]$.
\begin{algorithm}[H]
	\caption{Allocation for Non-Linear Costs}
	\label{optgalg}
	\begin{algorithmic}[1]
		\For {$i=n$ to $1$}
		\State $\hat{o}_i(x):=\sum_{j> i} \Big(o_j(\Delta_i-x+\sum_{k\mid i< k\leq j}(\Delta_k-x_k))- o_j(\sum_{k\mid i< k\leq j}(\Delta_k-x_k))\Big)$
		\State  $x_{i}=\argmin{x\geq 0}\, \max\{u_i(x),o_i(\Delta_i-x)+\hat{o_i}(x)\}$
		\EndFor
		\State {\bf Output:} $\{x_i\}_i$
	\end{algorithmic}
\end{algorithm}

\emph{Remark:} Given two sets of jobs $S_{k-1}=\{k-1,\dots,n\}$ and $S_k=\{k,\dots,n\}$. Let the output of Algorithm \ref{optgalg} over set $S_{k-1}$ be $\{x^{k-1}_i\}$, and over set $S_k$ be $\{x^{k}_i\}$. Then we have, $x^{k-1}_i=x^{k}_i$ for all $i\in S_k$. Also, note that the algorithm takes $n$ iterations, each involves a minimization of the maximum of two convex piece-wise linear functions (as $u_i(x)$ and $o_i(\Delta_i-x) + \hat{o}_i(x)$ are both piece-wise linear and convex in $x$). Therefore, each step involves finding the minimum of a convex piece-wise linear function. Hence, the algorithm runs in polynomial time (in $n$ and the number of pieces in the cost functions).  

\section{Robust Appointment Scheduling and Ordering}\label{sec:rass}

So far, we assumed that jobs are given to us in fixed order and our task was to find an optimal appointment schedule. Focusing on the case of constant per-unit costs, we now consider the joint problem of finding an ordering and an appointment schedule such that the resulting cost is minimized. More formally, consider a permutation $\pi$ over the set $[n]$, that determines the order of appointments. So given an ordering $\pi$, job $i$ is the $\pi(i)$-th appointment in the schedule. We let $A_i,C_i$ denote the start time and the completion time of the $i$-th appointment (or the $\pi^{-1}(i)$-th job). The joint scheduling and ordering problem can now be stated as,  
\begin{equation}
\min_{\pi:[n]\to [n], \{A_i\}_i}\quad \max_{t_i\in [0,\Delta_i]\, \forall i}\, \sum_{i=1}^n \max\{o_{\pi^{-1}(i)}(C_i-A_{i+1}),u_{\pi^{-1}(i)}(A_{i+1}-C_i)\}. \label{problemn}
\end{equation}
 Recall that under homogeneous underage costs $u_i=u\, \forall i$, we have a closed form solution for the problem for fixed permutation $\pi$. Letting $o(i)=\sum_{j\mid \pi(j)\geq \pi(i)} o_j$, the objective in this special case can be more simply stated as,
 \[ \min_{\pi:[n]\rightarrow [n]} \, \sum_{i=1}^n \frac{\Delta_io(i)u}{o(i)+u}.\]
 Let us call this problem RASO-H for brevity. Mittal et al.\ \cite{mittal} showed that RASO-H reduces to an instance of min-sum scheduling with concave costs $1||\sum w_jf(C_j)$. Here given jobs $j$ with processing time $p_j$, weight $w_j$, and a concave function $f(\cdot)$ over completion times $C_j$, the goal is to find an ordering that achieves the following,
 \begin{equation}\min_{\{C_j\}_j} \sum_{j=1}^n w_j f(C_j).\label{nonlin}
 \end{equation}
 Indeed, letting $w_j=\Delta_j$, $p_j=o_j$ and $f(C_j)=\frac{uC_j}{C_j+u}$, we see that any ordering $\pi$ for RASO-H is equivalent to an order $\pi'$ in, $1||\sum w_jf(C_j)$. Here $\pi'(j)=n-\pi(j)+1$ (i.e., orders are reversed as we move between the two problems). This reduces RASO-H to an instance of $1||\sum w_jf(C_j)$, while preserving the objective value. Therefore, algorithms for $1||\sum w_jf(C_j)$ can be used directly for RASO-H without any loss in guarantee. Using a result in \cite{smith}, the following rule,
 \[\text{ Smith's rule: Schedule jobs in the order of descending ratios $\frac{w_j}{p_{j}}$ (or ascending ratios $\frac{\Delta_j}{o_j}$),} \] is 1.137 approximate for RASO-H. In fact, there is an EPTAS for RASO-H due to the EPTAS for $1||\sum w_jf(C_j)$ \cite{eptas}. However, the hardness of RASO with a single underage cost, and more generally of $1||\sum w_jf(C_j)$, remains an intriguing open problem in scheduling theory \cite{eptas,megow,smith,mittal}.
 
 In the upcoming section, we consider the general problem for which no results were previously known. We give evidence indicating that the problem becomes much harder without the homogeneous underage costs assumption. %

\subsection{RASO with General Underage Costs}\label{rasog}
We start by showing that RASO is strongly NP-hard when there are at least two different underage costs via a reduction from the strongly NP hard problem of min-sum scheduling on identical parallel machines, $P||\sum w_jC_j$ (problem SS13 in Garey and Johnson \cite{garey}). Further, we also show (in the Appendix \ref{apxrasog}) that no list ordering rule (such as Smith's rule) can be better than $\Omega(n)$ approximate for the problem, ruling out the existence of ``simple'' approximation heuristics. We first state a property of the closed form allocation $x_i=\frac{\Delta_io(i)}{u_i+o(i)}$ that will be useful later on. 
 \begin{lem}\label{formulaworstcase}
 	The allocation given by $x_j=\frac{o(j)\Delta_j}{u_j+o(j)}\, \forall j\in[n]$ has all jobs taking maximum time and all taking minimum time as worst-cases.
 \end{lem}
 \begin{proof}
 	For the given allocation, the cost of all jobs taking zero time is the same as the cost of all jobs taking maximum time. We prove that these are worst-cases by induction on the number of jobs $n$. Claim is trivially true for a single job and suppose true for all $n\leq k-1$, then we prove for the case of $k$ jobs. Let $x'_i$ be the allocation given by the formula for the subset $[k]\backslash {1}$ of jobs. $x'_i=x_i$ for every $i\in[k]\backslash {1}$ and by assumption, all jobs underaged and all overaged are worst-cases for allocation $x'_i$. Therefore, if there exists a worst-case for allocation $x_i$ where job $1$ is underaged and takes zero time, then all jobs taking zero time is a worst-case and we're done.
 	Alternatively, if there is a worst-case with job $1$ overaged, we have from Lemma \ref{shift} that all jobs overaged is a worst-case.
 \end{proof}


 \begin{thm}\label{hard}
	RASO with underage costs $u_i\in \{v,w\}$ for all $i$, is strongly NP hard.
\end{thm}
\begin{proof}
	We reduce from min-sum scheduling on identical parallel machines . Consider an instance of the same with $m$ machines $j\in \mathcal{J}$, $n$ jobs $i\in \mathcal{P}$, weights $w_i$ and processing times $p_i$ (we assume integral values without affecting hardness). The objective is to find an assignment of jobs to machines that minimizes the weighted sum of completion times $\sum_{j\in\mathcal{J}} \sum_{i\in \mathcal{P}|i \text{ assigned to }j} w_iC_i$. Note that given an assignment, the job order on each machine is optimally determined by Smith's rule. Also let, $P=\sum_{i\in \mathcal{P}} p_i$ and $W=\sum_{i\in \mathcal{P}} w_i$.
	
	Now, we define an instance of RASO for $n+m-1$ jobs $i\in[n+m-1]$. The first $n$ jobs $i\in[n]$ will correspond to to jobs $\mathcal{P}$ and the remaining $m-1$ jobs will effectively model $m-1$ parallel machines. Let $u_i=M$ (large number to be decided later), $o_i=p_i$ and $\Delta_i=w_i$ for $i\in [n]$. For the remaining $m-1$ jobs, let $u_i=o_i=\Delta_i=0$. 
	Now, given an order $\pi$ (job $i$ is the $\pi(i)$th scheduled job) define $m_i=\pi^{-1}(\min \{j|j\geq \pi(i), \pi^{-1}(j)>n\})$ i.e., the first job $j$ out of jobs $[n+m-1]\backslash[n]$, that is scheduled on or after $i$ in ordering $\pi$ i.e., $\pi(j)\geq \pi(i)$. Note, $m_i$ could be ill defined for the last few jobs in the order. In any such case, set it to $\pi^{-1}(n+m)$  and let $\pi^{-1}(n+m)= n+m$. The idea is that jobs in $[n+m-1]\backslash[n]$ have zero underage cost and therefore, an optimal allocation will allocate enough time to every job $j\in[n+m-1]\backslash[n]$ so that the jobs succeeding $j$ are effectively isolated from jobs preceding $j$. In fact, we will see that every job $i\in[n]$ is in effect getting assigned to `machine' $m_i-n$. Further, let $o'(i)=\sum_{j|\pi(m_i)\geq \pi(j)\geq \pi(i)}o_j$ (note the difference from $o(i)=\sum_{j|\pi(j)\geq \pi(i)}o_j$). Let $S(\pi)=\sum_{i\in[n]}\Delta_io'(i)$. Observe that $S(\pi)$ is exactly the cost of assigning job $i$ to machine $m_i-n\in\{1,\dots,m\}$ in the parallel machine scheduling problem (with the order of the jobs on each machine given by reversing $\pi$), and every assignment of jobs to machines (with any order on each machine) is captured by some order $\pi$. Let $OPT(\pi)$ denote the worst-case cost of an optimal allocation given order $\pi$. We claim that for $M=P^2W+1$, we have $S(\pi)-1<OPT(\pi)\leq S(\pi)$. 
	
	For the upper bound, consider the allocation with $x_i=0$ for $i\in[n]$ and $x_i=\sum_{j|j\leq n, m_j=i}\Delta_j$ for $i\in[n+m-1]\backslash[n]$. Then, inductively in every time profile, every job $i\in[n+m-1]\backslash[n]$ ends by the end of its slot and the job following it starts at its appointed time. Therefore, job $i\in[n]$ is overaged by at most $\sum_{j|m_j=m_i,\pi(j)\leq \pi(i)}\Delta_j$ and has no underage cost. Summing up, the worst-case cost of this allocation is at most $\sum_{i\in[n]}o_i\sum_{j|m_j=m_i,\pi(j)\leq \pi(i)}\Delta_j=S(\pi)$. 
	
	For the lower bound, consider an arbitrary allocation and let $Y$ denote the total time allotted to jobs in $[n]$. The cost of all jobs taking zero time is at least $MY$. Now, consider the case of all jobs taking maximum time. Given job $i\in[n]$, let $s$ be the first job in the order with $m_s=m_i$. Only jobs $j\in[n]$ with $m_j=m_i$ are scheduled between $s$ and $i$, therefore, $C_{\pi(i)}\geq A_{\pi(s)}+\sum_{j|m_j=m_i,\pi(j)\leq \pi(i)}\Delta_j$ and $A_{\pi(i)+1}\leq A_{\pi(s)}+Y$. So the cost of job $i$ is at least $o_i(\sum_{j|m_j=m_i,\pi(j)\leq \pi(i)}\Delta_j-Y)$ and the total cost is at least $S(\pi)-PY$. Hence, $OPT(\pi)\geq \max\{MY,S(\pi)-PY\}\geq \frac{M}{M+P}S(\pi)>S(\pi)-1$, where the last inequality uses $S(\pi)\leq PW<M$. 
	
	Due to integrality of $w_i,p_i$, $S(\pi)$ is an integer and therefore, an optimal solution to this instance of RASO gives an optimal solution for the machine scheduling problem. Since $M$ is polynomially bounded in the input size (in unary), this proves strong NP hardness.  
\end{proof}
 
To develop approximation heuristics for RASO, unlike RASO-H, we cannot rely on a closed-form solution for the optimal allocation problem to simplify the problem. However, we do have a closed-form solution with cost guaranteed to be within a constant factor of the optimal due to Theorem \ref{closed2}. We further outline this idea and the obstacle we encounter in Appendix \ref{apxrasog2}. In a different approach, we consider the closed form allocation given by the formula $x_i=\frac{\Delta_i o(i)}{u_i+o(i)}$. Given order $\pi$, let $o(i)=\sum_{j\mid \pi(j)\geq \pi(i)} o_j$. For this allocation rule the ordering problem simplifies to,
\begin{equation}
\min_{\pi:[n]\rightarrow [n]} \sum_{j=1}^n \frac{\Delta_jo(j)u_j}{o(j)+u_j}.\label{genshed}
\end{equation}

This is equivalent to a min-sum scheduling problem of the form, 
\[\min_{\{C_j\}_j} \sum_j f_j(C_j).\]
Here $C_j$ is the completion time of job $j$ which has processing time $p_j=o_j$. Functions $f_j(C_j)=\Delta_ju_j \frac{C_j}{C_j+u_j}$ are concave, but now we have a different function for each job. 
We show (in Appendix \ref{apxrasog}) that given an $\alpha$ approximation for $1||\sum_jf_j(C_j)$, there is a $4\alpha n$ approximation for RASO. The additional factor of $4n$ arises from the fact that (\ref{genshed}) does not represent the true objective \eqref{problemn}. Recall that the closed form allocation formula $x_j=\Delta_j \frac{o(j)}{o(j)+u_j}$, is not only suboptimal but can be arbitrarily worse than the optimal allocation for certain orders. So it is perhaps surprising that using this allocation we can still achieve some approximation bound. Note that for the min-sum scheduling problem, $1||\sum f_j(C_j)$,  \cite{mestre} gives a $4+\epsilon$ approximation for arbitrary $f_j$. 
Using this algorithm, we have a $(16+\epsilon)n$ approximation for RASO. 
We also show (in Appendix \ref{apxrasog}) that no heuristic that orders jobs based on a simple ratio can beat $\Omega(n)$. Here we use the term `simple ratio' to refer to any real valued function evaluated independently for each job, using only the job parameters ($o_i,\Delta_i,u_i$). A list ordering heuristic for this function then simply orders jobs in ascending or descending order of the function values.

\subsection{Homogeneous Underage Costs}\label{rasohomo}
 In this section, we develop some easy to implement heuristics with improved approximations for RASO-H. As we mentioned earlier, the computational complexity of RASO-H and more generally, $1||\sum w_jf(C_j)$ with concave cost function $f$, remains open but, there is a scaling based EPTAS for RASO-H due to the EPTAS for $1||\sum w_jf(C_j)$. The EPTAS is non-trivial and not easy to implement in practice. While heuristics such as Smith's rule offer a more practical alternative, they do not necessarily use any of the structure specific to RASO-H. For example, Smith's rule is suboptimal for RASO-H even with just two jobs. 
 
 This motivates us to consider the following heuristic, 
\[
\text{ Customized-Smith's rule (C-Smith): Schedule jobs in ascending order of $\frac{\Delta_i}{o_i(o_i+u)}$. }\]
 This straightforward heuristic is optimal for two jobs by design, and has an approximation guarantee of $\beta$, where $1.06036<\beta<1.06043$ i.e., C-Smith order costs at most $1.06043$ times the optimal cost, and there exists an instance where it costs at least $1.06036$ times the optimal cost. 
While the order output by C-Smith is optimal for two jobs, we could more generally seek an order that is optimal w.r.t.\ exchanging the order of any two consecutive jobs. Let us call such an order \emph{locally optimal}. The schedule output by both Smith and C-Smith is not locally optimal. In fact, as Lemma \ref{pairorder} below shows, whether two adjacent jobs are in the right relative order depends on the jobs scheduled after them, so a list-ordering heuristic need not be locally optimal. 
The algorithm we discuss next outputs a locally optimal order (here and in what follows, we let $u=1$ w.l.o.g., by scaling the overage costs). While not as simple as C-Smith, it is still fairly easy to implement -- at each step, it computes a new set of ratios for the remaining jobs and picks the one with the best ratio, removing it from further consideration. However, a naive implementation has runtime $O(n^2)$, in contrast to $O(n\log n)$ for list ordering heuristics.

\begin{algorithm}[H]
	\caption{Locally Optimal Algorithm}
	\label{localalg}
	\begin{algorithmic}[1]
		\State {\bf Initialize:} $S=\emptyset,o(S)=0$
		\For {$i=1$ to $n$}
		\State Find, \[j=\argmax{k\in[n]\backslash S} \frac{\Delta_k}{o_k(o_k+o(S)+1)}.\] \[\text{In case of a tie, pick the job with largest overage cost.}\]
		\State $\pi(j)=n-i+1,S=S\cup \{j\}, o(S)=o(S)+o_j$
		\EndFor
		\State {\bf Output:} $\pi(.)$
	\end{algorithmic}
\end{algorithm}
We will show that the above algorithm outputs a solution that is at least as good as C-Smith on every instance. Thus, it is also $\beta$ approximate. There exist instances for which the algorithm is not optimal (see Appendix \ref{apxrasohomo}). However, we leave finding the exact guarantee of the algorithm as an open problem.

Let us now understand why Algorithm \ref{localalg} outputs a locally optimal order. For this we need the following lemma which is proven in Appendix \ref{apxrasohomo}.
\begin{lem}\label{pairorder}
	Given an arbitrary schedule $\pi$, with two jobs $i,j$ scheduled consecutively i.e., $\pi(j)=\pi(i)+1$. Let $o_x=\sum_{\pi(k)> \pi(j)} o_k$. Then the order is optimal w.r.t. swapping jobs $i,j$ iff $$\frac{\Delta_j}{o_j(o_j+o_x+1)}\geq \frac{\Delta_i}{o_i(o_i+o_x+1)}.$$ 
\end{lem}

Observe that when we just have two jobs, $o_x=0$ and the optimal order is $i$ followed by $j$ iff $\frac{\Delta_j}{o_j(o_j+1)}\geq\frac{\Delta_i}{o_i(o_i+1)}$, which is the C-Smith order. Also note that if $o_x\to \infty$, then we have $i$ before $j$ iff $\frac{\Delta_j}{o_j}\geq \frac{\Delta_i}{o_i}$, indicating that Smith's rule can be optimal for ordering jobs earlier in the schedule if the overage costs are large enough.

\begin{lem}
	Algorithm \ref{localalg} is locally optimal.
\end{lem}
\begin{proof}
	Consider arbitrary jobs $i$ and $j$ with $\pi(j)=\pi(i)+1$ in the order $\pi$ output by the algorithm. Let $o_{S_j}$ represent the value of $o(S)$ at the iteration where job $j$ was picked out. We have, $\frac{\Delta_j}{o_j(o_j+o_{S_j}+1)}\geq \frac{\Delta_{i}}{o_{i}(o_{i}+o_{S_j}+1)}$. Lemma \ref{pairorder} now implies that the current order is locally optimal w.r.t.\ swaps between any two consecutive jobs in order $\pi$.
\end{proof}

Next, we examine how the optimal ordering of a pair of jobs changes due to other preceding jobs in a given schedule. This will allow us to prove that the order output by Algorithm \ref{localalg} is always at least as good as the order given by C-Smith. The proof of the following lemma is included in Appendix \ref{apxrasohomo}.
\begin{lem}\label{orderimpl}
	If $\frac{\Delta_i}{o_i(o_i+1)}\geq\frac{\Delta_j}{o_j(o_j+1)}$ and $\frac{\Delta_j}{o_j(o_j+x+1)}\geq\frac{\Delta_i}{o_i(o_i+x+1)}$ for some $x>0$, then we have, $o_j\geq o_i$. In a similar vein, given $o_j\geq o_i$ and $\frac{\Delta_j}{o_j(o_j+x+1)}\geq\frac{\Delta_i}{o_i(o_i+x+1)}$, we have, $\frac{\Delta_j}{o_j(o_j+x'+1)}\geq \frac{\Delta_i}{o_i(o_i+x'+1)}$ for every $x'\geq x$.
\end{lem}

\begin{coro}\label{ordercoro}
	Suppose $\frac{\Delta_i}{o_i(o_i+1)}=\frac{\Delta_j}{o_j(o_j+1)}$. If $o_j\geq o_i$, then $\frac{\Delta_j}{o_j(o_j+x+1)}\geq \frac{\Delta_i}{o_i(o_i+x+1)}$ for every $x\geq 0$. Conversely, if $\frac{\Delta_j}{o_j(o_j+x+1)}\geq \frac{\Delta_i}{o_i(o_i+x+1)}$ for some $x>0$ then $o_j\geq o_i$.
\end{coro}
\begin{thm}
	 For every instance of RASO-H, Algorithm \ref{localalg} outputs an order that is at least as good as the order given by C-Smith. 
\end{thm}  	
\begin{proof}
	Let $\pi_1$ be the order from C-Smith and $\pi_2$ the order from Algorithm \ref{localalg}. Further, suppose that the orders coincide on the last $i$ jobs in the schedules. So $n-i$ is the largest index such that $\pi^{-1}_1(n-i)\neq \pi^{-1}_2(n-i)$. Let $k=\pi^{-1}_2(n-i)$, then clearly, $\pi_1(k)<n-i$. We will show that the order which schedules $k$ at $n-i$, but otherwise follows the C-Smith order $\pi_1$, is at least as good as $\pi_1$ and more generally, we demonstrate a way to transform $\pi_1$ to $\pi_2$ while never increasing the cost.
	
	We have $\frac{\Delta_{j}}{o_{j}(o_{j}+1)}\geq\frac{\Delta_k}{o_k(o_k+1)}$ for every job $j$ scheduled after $k$ in order $\pi_1$ i.e., $\pi_1(j)\geq \pi_1(k)$ (due to C-Smith). Let $o(n-i+1)=\sum_{l\mid \pi_1(l)\geq n-i+1} o_l=\sum_{l\mid \pi_2(l)\geq n-i+1} o_l$, and let $J=\{j\mid \pi_1(k)\leq \pi_1(j)\leq n-i\}$. Since $\pi_2(k)=n-i$, we have that for all jobs $j\in J$, \[\frac{\Delta_k}{o_k(o_k+o(n-i+1)+1)}\geq\frac{\Delta_j}{o_j(o_j+o(n-i+1)+1)}.\] It follows that if $o(n-i+1)=0$, then due to the earlier inequality $\frac{\Delta_{j}}{o_{j}(o_{j}+1)}\geq\frac{\Delta_k}{o_k(o_k+1)}$, we in fact get, $\frac{\Delta_k}{o_k(o_k+1)}=\frac{\Delta_j}{o_j(o_j+1)}$ for  all $j\in J$. Then, from the tie breaking rule in Algorithm \ref{localalg} we have, $o_k\geq o_j$ for all $j\in J$. On the other hand, if $o(n-i+1)>0$ then from Lemma \ref{orderimpl} we again have $o_k\geq o_j$. So using Corollary \ref{ordercoro} and Lemma \ref{orderimpl} (if $o(n-i+1)$ is non-zero) we have,  \[\frac{\Delta_k}{o_k(o_k+\sum_{l\mid \pi_1(l)> \pi_1(j)} o_l+1)}\geq\frac{\Delta_j}{o_j(o_j+\sum_{l\mid \pi_1(l)> \pi_1(j)} o_l+1)},\] for every $j\in J$. Therefore, when $j=\pi^{-1}_1(\pi_1(k)+1)$, we have that swapping $j$ and $k$ in $\pi_1$ does not result in a worse schedule. Subsequently, we can swap $k$ with every $j\in J$ in $\pi_1$ to get a new schedule $\pi'_1$, where $\pi'_1(k)=n-i=\pi_2(k)$. The order in $\pi'_1$ is otherwise the same as $\pi_1$. The cost of $\pi'_1$ is at most the cost of $\pi_1$, and the first index where $\pi'_1$ and $\pi_2$ differ is $n-i-1$. Repeating this process a finite number of times proves the main claim. 
\end{proof}

It remains to show the $\beta$ approximation for C-Smith.
For this analysis it will be much more convenient both for clarity as well as notation, to focus on analyzing C-Smith for the scheduling problem $1||\sum_j w_jf(C_j)$, where $f(C_j)=\frac{C_ju}{C_j+u}$ and $w_j=\Delta_j$,$p_j=o_j$. Note that we can let $u=1$ w.l.o.g. Also recall, there is a cost preserving bijection between orders for $1||\sum_j w_jf(C_j)$ and RASO-H -- reversing the order when moving from one to the other. Therefore, we will show the following equivalent theorem.
\begin{thm}\label{maint}
    For some constant $\beta\in [1.06036,1.06043]$, scheduling jobs in descending order of $\frac{w_j}{p_j(p_j+1)}$ is exactly $\beta$ approximate for the scheduling problem $1||\sum_j w_j\frac{C_j}{C_j+1}$. 
\end{thm}
To prove the above theorem, we establish several intermediate results that characterize and simplify the worst-case instance for C-Smith. First, we show in Lemma \ref{normalize} that there is a worst-case instance for C-Smith where all jobs are tied and in fact have ratio 1. This is shown via a generalization of Lemma~3.5 in \cite{smith} (Lemma 21 in \cite{eptas}). Then, in Lemma \ref{order} we show that in the worst-case C-Smith orders jobs in ascending order of processing times and the optimal order is the exact reverse. Interestingly, this is the opposite of Lemma~3.6 in \cite{smith} (Proposition~20 in \cite{eptas}) for Smith's rule, where the optimal order is ascending in $p_i$ and the worst order is descending. Given these properties, we can formulate a non-convex optimization problem in infinitely many variables, the optimal value of which is the approximation ratio, and every optimal solution is a worst-case instance. To then get lower and upper bounds on the optimum value, in Lemmas \ref{lbp} to \ref{sumbound} we utilize properties specific to our objective $f(x)=\frac{x}{x+1}$ to approximate the problem (which has infinitely many variables) with a family of optimization problems, that while still non-convex, have a finite number of variables. More complex objectives from the family give a tighter upper bound on true approximation ratio, but the number of variables involved increases. We find the global optimum to a problem in the family with five variables, and this closely matches our lower bound. We solve such a non-convex problem to global optimality by establishing upper and lower bounds on variables and using linear cuts, both of which allow us to then effectively use a nonlinear globally optimal MINLP solver, Couenne \cite{cou}. The details are provided in Appendix \ref{apxcsmith}.

\section{Conclusion \& Open Problems}\label{conclusion}
We considered the robust appointment scheduling problem with general underage costs. For the appointment scheduling problem with fixed jobs order, we found a simple LP that gives a 2-approximation for the problem under constant per-unit costs. Then we further refined this LP, resulting in a closed form solution for optimal allocations in special cases (generalizing previous results in the robust model, and complementing similar results in other models). We also showed that our results and approach extend more generally to convex piece-wise linear costs. When seeking an optimal allocation for the general case using our approach, more complications arise and it is not clear if one can still construct a linear (or convex) program such that an optimal solution to the program is also an optimal allocation. We leave finding an optimal solution for the general case (or showing hardness), as an open problem. 

%
In the second setting, we considered the problem of jointly finding the optimal order and allocation given that order for the case of constant per-unit costs. For the case of heterogeneous underage costs, we show that the problem is strongly NP hard and no list ordering policy can do better than $\Omega(n)$ to approximate the optimal value. We also gave a heuristic that achieves this bound. Finding a better approximation for this setting remains another interesting open problem. For the case of homogeneous underage costs, we designed two simple and practical heuristics that are guaranteed to be within $\approx 1.06$ of the optimal.  
\section*{Acknowledgements}
The authors would like to thank James B. Orlin for helpful discussions. This work was partially supported by ONR Grant N00014-17-1-2194, NSF Grant CMMI 1636046 and Alexander von Humboldt Foundation (with funds from BMBF).

\bibliographystyle{plain}
\bibliography{rapt}
\begin{appendix}
\section{2-approximation for RAS with Non-decreasing, Convex Piece-wise Linear Costs}\label{ext1}
Given non-decreasing costs $o_i(\cdot),u_i(\cdot)$, we have the following lemmas. Note that we do not need the additional assumptions of piece-wise linearity and convexity for the below.
\begin{lem}\label{sanity1g}
	 There exists an optimal allocation $\{x_i\}_i$, such that,  $\sum_{j=1}^{k}x_{j}\leq \sum_{j=1}^{k} \Delta_j$\ for all $k = 1, \ldots, n$. 
\end{lem}
\begin{proof}
	Given an optimal allocation $\{x_i\}_i$, let $\delta_k=\sum_{j=1}^{k}x_{j}- \sum_{j=1}^{k} \Delta_j$, for $k\in[n]$. Assume that $\{x_i\}_i$ is an optimal allocation such that, there is no optimal allocation $\{y_i\}_i$ with $\sum_i y_i<\sum_i x_i$. Suppose now, $\delta_k>0$ for some $k$, and let $k_0$ be the smallest such $k$. Given this, notice that even if the first $k_0$ jobs take maximum time, job $k_0$ can never be overaged. Therefore, $k_0$ is underaged in every worst-case and decreasing $x_{k_0}$ by a small amount gives another optimal allocation. In particular, the allocation $\{x'_i\}_i$ where, $x'_{k_0}=x_{k_0}-\delta_{k_0}$ and $x'_i=x_i$ for all other $i$, is optimal. This is a contradiction since $\sum_i x'_i<\sum_i x_i$. 
\end{proof}
\begin{lem}\label{costubg}
	Given an allocation $\{x_i\}_{i}$ where $\sum_{j=1}^{i}x_{j}\leq \sum_{j=1}^{i} \Delta_j$\ for all $i$, and a time profile $T=(t_i)$. The cost $ c(T)$ of profile $T$ is at most $\sum_{i=1}^n o_i(\sum_{j\leq i}(\Delta_j-x_j))+\sum_{i=1}^n u_i (x_i)$.  
\end{lem}
\begin{proof}
	The result follows by observing that for any given time profile $T$, job $i$ is either underaged with $c(T,i)\leq u_i(x_i)$ or overaged with cost at most $o_i(\sum_{j \leq i} (\Delta_j-x_j))$. The latter follows from the fact that $C_i\leq \sum_{j \leq i}\Delta_j$ since $\sum_{j=1}^{k}x_{j}\leq \sum_{j=1}^{k} \Delta_j$ for all $k\in[i]$. 
\end{proof}

Using the above results, the following program now gives us a 2-approximation for the general problem. 
\begin{eqnarray}
  &\min & \sum_{i=1}^n o_i(\sum_{j\leq i}(\Delta_j-y_j))+\sum_{i=1}^n u_i (y_i) \nonumber\\
&{\rm s.t.} 
 & \sum_{j=1}^{k}(y_j-\Delta_j)\leq 0\quad \forall k\in [n];\quad \label{auxg}\\
&& y_j\geq 0 \quad \forall j\in [n].\nonumber
\end{eqnarray} 

\begin{thm}\label{closed2pie}
	The optimal solution to the above program is a (tight) $2$-approximation to the optimal allocation. 
\end{thm}
\begin{proof}
 From constraints~\eqref{auxg} and Lemma \ref{costubg}, the objective is an upper bound on the worst-case cost of the allocation given by any feasible solution to the program. Let $\{y^*_i\}$ denote an optimal solution to the program. Lemma \ref{sanity1g} implies that every optimal allocation $\{x_i\}$ is a feasible solution to the program. Now, the worst-case cost of an optimal allocation, denoted $OPT$, is at least $\max\{\sum_iu_i(x_i),\sum_i o_i(\sum_{j\leq i}(\Delta_j-x_j))\}$ (we show later that these two costs are both in fact, equal to $OPT$). 
	Therefore, \[ \sum_iu_i(y^*_i) +\sum_i o_i(\sum_{j\leq i}(\Delta_j-y^*_j))\leq \sum_iu_i(x_i)+\sum_i o_i(\sum_{j\leq i}(\Delta_j-x_j))\leq 2~OPT.\]
	For a tight instance, refer to the proof of Theorem \ref{2apx}. If the costs $o_i(\cdot),u_i(\cdot)$ are piece-wise linear and convex, the program can be solved efficiently using polynomial time algorithms for network flows with convex piece-wise linear arc costs, such as we argued in Section \ref{2apxg}.
\end{proof}

\section{Missing Proofs from Section \ref{sec:optspec}}\label{apxoptspec}
\begin{proof}[Proof of Lemma \ref{l1}]
	Let us fix optimal allocation $\{x_i\}_i$ and arbitrary job $i$. We prove the claim in parts by considering three cases; when $i$ is \emph{strictly} overaged or \emph{strictly} underaged in all worst-cases, when it is underaged in all cases, and finally when it is overaged in all cases. Recall that being charged for overage means job is completed after time $A_{i+1}$, which could occur even if $t_i=0$.
	
	First, consider the case where $i$ is strictly underaged in all worst-cases. Then reducing $x_i$ leads to a better allocation, yielding a contradiction. Note that $x_i>0$ since we assumed $i$ is strictly underaged. When $i$ is always strictly overaged, let $k$ be the smallest index larger than $i$ such that, $x_k>0$. We can reduce the cost in every worst-case by increasing $x_i$ by some small $\epsilon>0$ and decreasing $x_k$ by $\epsilon$. Note that if we only increase $x_i$ we might underage some job $j>i$. Decreasing $x_k$ by the same $\epsilon$ prevents this from occurring. Therefore, in an optimal allocation no job can be strictly overaged (underaged) in every worst-case.  
	
	Now, consider the scenario where job $i$ is always underaged but there is a worst-case $T=(t_j)$ where $i$ has zero cost, so $t_i=0$ and $C_i=A_{i+1}$. In fact, let $i$ be the job with largest index that has this property. We also have $\Delta_i>0$, since if $\Delta_j=0$ for some job $j$ and there is a worst-case where $j$ is underaged with zero cost, it is also overaged (by definition). Suppose for every job $k$ with $k>i$, $x_k=0$ (or $i=n$). Then, since $\Delta_i>0$ we have that profile $T'$ with $t'_i=\Delta_i$ and $t'_j=t_j$ for $j\neq i$, is also a worst-case. Therefore, assume that $x_k>0$ for some $k>i$ and let $k$ denote the smallest such index. Now, there exists a worst-case $Z=\{z_j\}_j$ where $k$ is overaged, since no job is strictly underaged in all worst-cases, and $i$ is the largest index job which is underaged in all worst-cases. Given such a profile $Z$ where $k$ is overaged, consider the subset of jobs given by $S_{i+1}=\{i+1,\dots,n\}$. Since $i$ is underaged in $Z$ by assumption, we have that that the time profile represented by $\{z_j\}_{j\in S_{i+1}}$ is a worst-case for allocation $\{x_j\}_{j\in S_{i+1}}$. Therefore, it follows that the profile $Q=\{q_j\}_{j\in [n]}$ where $q_j=z_j$ for $j\geq i+1$ and $q_j=t_j$ otherwise, is a worst-case for allocation $\{x_j\}_{j\in[n]}$. Modifying $q_i=t_i=0$ to $q'_i=\epsilon>0$ and reducing $q_j>0$ to $q'_j=q_j-\epsilon$ for the smallest $i<j\leq k$ (such a $j$ exists since $k$ is overaged in $Q$ and $x_k>0$), yields a worst-case $Q'$ with $i$ overaged, contradiction. 
	
	Finally, consider the remaining case where job $i$ is always overaged but there is a worst-case where $i$ is overaged with zero cost. Since $i$ completes strictly after $A_{i+1}$ in this worst-case (the case of zero cost where $C_i=A_{i+1}$ counts as underage), this case only occurs if $o_i=0$, which we excluded by assumption.
	
	We end with the case where $o_j=0$ for some jobs $j$. Replace each such $o_j$ by $\epsilon>0$ and let $\{x^\epsilon_i\}_i$ denote an optimal allocation for the resulting instance. Due to Lemma \ref{sanity1}, we may assume $x^\epsilon_i\in[0,\sum_j\Delta_j]$ for every $i$. The worst-case cost of an allocation is the maximum of a continuous function over a fixed compact set of time profiles, and is therefore continuous in the allocation and the cost parameters. Thus, any limit point $\{x^*_i\}_i$ of $\{x^\epsilon_i\}_i$ as $\epsilon\to0$ is an optimal allocation for the original instance. Further, from Lemmas \ref{worst2} and \ref{worst1} the worst-case cost of $\{x^\epsilon_i\}_i$ equals both $\sum_i u_ix^\epsilon_i$ and $\sum_i o(i)(\Delta_i-x^\epsilon_i)$ (the latter is the cost of all jobs overaged due to Lemma \ref{sanity1}). Taking limits, all jobs underaged and all jobs overaged are worst-cases for $\{x^*_i\}_i$. Therefore, Lemmas \ref{worst2} and \ref{worst1} hold for at least one optimal allocation even if $o_j=0$ for some $j$.
\end{proof}

Now we state and prove two lemmas that are then used to prove Lemma \ref{l3}; one of the key Lemmas for formulating the strengthened LP-2.
\begin{lem}\label{monotone}
	Recall $x^S_i$ denotes the time allocation for job $i$ in an optimal allocation over jobs in set $S$ and $S_i=\{i,\dots,n\}$ (here, we consider optimal allocations for which all jobs underaged and all jobs overaged are worst-cases). We have, $\sum_{i\in S_k} u_ix^{S_k}_i \leq \sum_{i\in S_k}u_i x^{[n]}_i$ for every $k\in[n]$.
\end{lem}
\begin{proof}
Note that the LHS denotes the worst-case cost of an optimal allocation over $S_{k}$. Now, consider the possibly sub-optimal allocation $x^{[n]}_i$ for jobs $S_{k}$. The RHS denotes the cost of all jobs underaged under this allocation, however a priori this need not be a worst-case cost for the allocation. So to show the inequality we prove that all jobs underaged is a worst-case for the allocation in RHS.
	
	Assume this is not the case i.e., for allocation $x^{[n]}_i$ over $S_k$ there is a profile, denoted by $Y$, where some job takes non-zero time, and $c(Y)>\sum_i u_i x^{[n]}_i$. Now consider the set of jobs $[n]$ and the now optimal allocation $x^{[n]}_i$ over this set. We have that $\sum_{i\in [n]} u_ix^{[n]}_i< \sum_{i<k}u_i x^{[n]}_i +c(Y)$. This implies that for allocation  $x^{[n]}_i$ over $[n]$, the profile where jobs $1$ to $k-1$ take zero time and other jobs are as in $Y$ is strictly worse than underaging all jobs, contradicting Lemma \ref{worst2}.
\end{proof}

\begin{lem}\label{costlb}
	Given an arbitrary allocation $\{x_i\}_{i}$ and a time profile $T=(t_i)$ for $i\in[n]$, the cost $c(T,i)$ of job $i$ is at least $o_i\sum_{j \leq i} (t_j-x_j)$. 
\end{lem}
\begin{proof}
	We have that $c(T,i)\geq  o_i(C_i-A_{i+1})$ by definition. Also, $ C_i\geq \sum_{j\leq i} t_j$ and $A_{i+1}=\sum_{j\leq i} x_j$. Therefore, $c(T,i)\geq o_i \sum_{j\leq i}(t_j-x_j)$.
\end{proof}

Now, we have the ingredients required to prove Lemma \ref{l3}.
\begin{proof}[Proof of Lemma \ref{l3}]
Observe that for $k=1$ we have equality from~\eqref{sumcost}. Hence, it suffices to prove the equivalent inequalities, \[\sum_{j=1}^{k} (u_j+o(j))y_j\leq \sum_{j=1}^{k} o(j)\Delta_j\quad  \text{for all } k\in[n].\] Our proof is via induction on $k$, with $n>1$ fixed. For the base case of $k=1$, suppose the statement is false i.e., $u_1x_1>o(1)(\Delta_1-x_1)$. We derive a contradiction by constructing a new allocation $\{x'_i\}_i$ that is better than the optimal $\{x_i\}_i$. Let, 
	\begin{equation}
	x'_1=\frac{o(1)\Delta_1}{u_1+o(1)}< x_1, \label{defxn}
	\end{equation}
	and $x'_i=x^{S_2}_i$ for every $i\geq 2$ (recall $\{x_i^{S_2}\}$ denotes an optimal allocation over jobs in subset $S_2=\{2,\dots,n\}$). By design, in $\{x'_i\}_{i\in[n]}$ the total cost when all jobs take maximum time equals total cost when all jobs take zero time,
	\[	\sum_{i=1}^n o(i)(\Delta_i-x'_i) =o(1)(\Delta_1-x'_1)+\sum_{i=2}^{n} o(i)(\Delta_i-x'_i)
	=u_1 x'_1 + \sum_{i=2}^{n} u_ix'_i.\]
	Here the second equality follows from~\eqref{sumcost} for optimal allocations over $S_2$ and (\ref{defxn}). Moreover, this cost is strictly smaller than worst-case cost in allocation $\{x_i\}_i$ i.e., $\sum_{i=1}^n u_ix'_i < u_1x_1+\sum_{i=2}^{n} u_i x'_i\leq \sum_{i=1}^n u_i x_i$, where the last inequality follows from Lemma \ref{monotone}. To finish the proof for the base case we need to show that all jobs taking maximum time or alternatively all taking zero time is a worst-case in the new allocation. We argue this by considering two cases based on job $1$, while using the fact that $\{x'_i\}_{i\in S_2}$ is optimal for the set of jobs $S_2$. First, if there is a worst-case where $1$ is overaged, then all jobs taking maximum time is a worst-case due to Lemma \ref{shift}. Similarly, if job $1$ is underaged in some worst-case then all jobs underaged is a worst-case due to Lemma \ref{worst2}. This completes the base case. 
	
	Next, assume we are given $\sum_{j=1}^{i} u_jx_j\leq \sum_{j=1}^{i} o(j)(\Delta_j-x_j)$ for every $i\leq k$ and we wish to show the inequality for $i=k+1$. 
	As before, we proceed by contradiction, so assume $\sum_{j=1}^{k+1} u_jx_j> \sum_{j=1}^{k+1} o(j)(\Delta_j-x_j)$. Similar to the base case, we construct a new allocation $\{x'_i\}_i$ where, $$x'_{k+1}=\frac{1}{u_{k+1}+o(k+1)}[o(k+1)\Delta_{k+1}+\sum_{j=1}^{k} o(j)(\Delta_j-x_j)-\sum_{j=1}^{k} u_jx_j]< x_{k+1}.$$
	Further, $x'_i=x_i$ for every $i\leq k$ and $x'_{i}=x^{S_{k+2}}_i$ for every job in $S_{k+2}$. Note that $x'_{k+1}\geq 0$ by the induction hypothesis, and $\{x'_i\}_i$ satisfies the conditions in Lemma \ref{sanity1} since $x'_{k+1}<x_{k+1}$. Moreover, similar to the base case, $\sum_i u_ix'_i<\sum_i u_ix_i$ since $x'_{k+1}<x_{k+1}$ and $\sum_{i\in S_{k+2}}u_ix^{S_{k+2}}_i\leq \sum_{i\in S_{k+2}}u_ix_i$ (Lemma \ref{monotone}). Same as the base case, observe that in the new allocation the cost when all jobs are underaged equals the cost when all jobs are overaged. We show that at least one of them is a worst-case (and thus both are). 
	
	If we have a worst-case for $\{x'_i\}_i$ with $k+1$ underaged, then there exists a worst-case with all jobs in $S_{k+1}$ underaged since Lemma \ref{worst2} implies underaging all jobs in $S_{k+2}$ is a worst-case for $\{x'_i\}_{i\in S_{k+2}}$. Let $Y$ be a worst-case for $\{x'_i\}_{i\in [n]}$ where jobs in $S_{k+1}$ are all underaged and assume there exists some job in $[k]$ that is overaged in $Y$ (otherwise we are done). Then we have that sum of costs from jobs in $[k+1]$ is strictly larger than the cost of underaging those jobs i.e., $\sum_{i=1}^{k+1} c(Y,i)> \sum_{i=1}^{k+1} u_ix'_i=\sum_{i=1}^{k} u_ix_i+u_{k+1}x'_{k+1}$. Therefore, $\sum_{i=1}^{k+1} c(Y,i)+u_{k+1}(x_{k+1}-x'_{k+1})>\sum_{i=1}^{k+1} u_ix_i$. Observe that this yields a contradiction for the allocation $\{x_i\}_{i\in [n]}$, since the cost due to all jobs underaged is strictly smaller than the cost when jobs in $[k]$ are as in $Y$ and jobs in $S_{k+1}$ are all underaged.
	
	Finally, suppose there is no worst-case for $\{x'_i\}_{i\in[n]}$ where $k+1$ is underaged. In other words, $k+1$ is overaged in all worst-cases. Then due to Lemma \ref{shift} and Lemma \ref{worst1}, there is a worst-case where all jobs in $S_{k+1}$ are overaged and take maximum time. Now, if there also exists a worst-case where all jobs in $[k]$ are overaged, then clearly all jobs overaged is a worst-case. So assume that in every worst-case, some job in $[k]$ is underaged. Then let us consider a worst-case, call it $Z$, where all jobs in $S_{j+1}$ are overaged, where $j$ is the largest index job in $[k]$ that is underaged in $Z$. W.l.o.g., jobs in $S_{j+1}$ take maximum time in $Z$, since increasing the time taken by an overaged job keeps it and all subsequent jobs overaged, and does not decrease the cost. Since $Z$ is strictly worse than the case of all jobs overaged, we have, $\sum_{i=1}^n o(i)(\Delta_i-x'_i)< \sum_{i=1}^n c(Z,i)$. Therefore,
	\begin{eqnarray}
	\sum_{i=1}^j o(i)(\Delta_i-x_i) <\sum_{i=1}^j c(Z,i). \label{partsum}
	\end{eqnarray} 
	The inequality above follows from $\sum_{i=j+1}^{n} c(Z,i)= \sum_{i=j+1}^{n} o(i)(\Delta_i-x'_i)$ and using the fact that $x'_i=x_i$ for $i\leq j\leq k$. Now, consider again the allocation $\{x_i\}_{i\in [n]}$ and a time profile $T$ where jobs in $[j]$ are as in $Z$ and jobs in $S_{j+1}$ all take maximum time. Note that jobs in $S_{j+1}$ need not all be overaged in this case. However, since job $j$ ends by time $A_{j+1}$ in $Z$ and $x'_i=x_i$ for $i\leq j$, job $j+1$ starts at $A_{j+1}$ in $T$. So for every $i\geq j+1$, $C_i\geq A_{j+1}+\sum_{l=j+1}^{i}\Delta_l$ and the cost of job $i$ is at least $o_i\sum_{l=j+1}^{i}(\Delta_l-x_l)$ (similar to Lemma \ref{costlb}). Thus, the net cost due to jobs in $S_{j+1}$ is at least $\sum_{i=j+1}^{n} o(i)(\Delta_i-x_i)$. Combining this with (\ref{partsum}), we have that for allocation $\{x_i\}_{i\in[n]}$, case $T$ is strictly worse than overaging all jobs, contradicting Lemma \ref{worst1}.

\end{proof}

With this, we are now ready to prove the result on closed form optimal allocations.
\begin{proof}[Proof of Theorem \ref{sidedspec}]
	We first show that the allocation given by the formula is a feasible solution to the LP and that all jobs underaged/overaged are worst-cases. Then we prove that this allocation is an optimal solution for the LP. 
	
	To check feasibility, first note that equation (\ref{bala}) is satisfied by definition, since each term $o(i)\Delta_i$ appears exactly once in the sum $\sum_j (o(j)+u_j)x_j$. Constraints (\ref{knap}) are also similarly easy to verify. However, $x_i$ can be larger than $\Delta_i$ for some $i$, so checking (\ref{aux}) is less trivial. Indeed, this is where we use the assumption $u_{m_i}\geq o(i)-o(m_i)$ for every $i$. Note that if for job $i$, $m_i\neq i$, then $x_i=0$ and $x_{m_i}$ has a term $\frac{\Delta_io(i)}{u_{m_i}+o(m_i)}\leq \Delta_i$ in the summation. More generally, we have that for every $k$ either $x_k=0$ or $x_k\leq \sum_{i\leq k|k=m_i}\Delta_i$. It now follows that constraints (\ref{aux}) are satisfied. 
	
	To show that all jobs underaged/overaged (both have same cost due to (\ref{bala})) is a worst-case, we proceed by induction on the number of jobs $n$. For $n=1$ this is trivial. Assuming the claim holds for $k-1$ jobs, we prove the statement for $k$ jobs. Let the allocation given by the formula for $k-1$ jobs given by $S_2=\{2,\dots,k\}$, be $\{x'_i\}_{i\in S_2}$. Let $\{x_i\}_{i\in [k]}$ be the allocation for jobs $i\in[k]$ and note that the values $m_i$ for $i\in S_2$ are the same in both cases. We will break the proof down into two cases. First, suppose $m_1=1$. Then $x_i=x'_i$ for every job $i\in S_2$. Therefore, if for allocation $x_i$ there exists a worst-case where job $1$ is underaged, then all jobs underaged is a worst-case (because all jobs underaged is a worst-case for allocation $x'_i$ by assumption). Otherwise, job $1$ is overaged in some worst-case and using Lemma \ref{shift} and the induction assumption, we have that all jobs overaged is a worst-case.
	
	Now, suppose $m:=m_1\neq 1$. Then $x_i=0$ for every $i< m$, $x_i=x'_i$ for all $i\in S_2\backslash \{m\}$ and $x_{m}=x'_{m}+\frac{\Delta_1o(1)}{u_{m}+o(m)}$. Observe that $\{x_i\}_{i\in S_2}$ is exactly the allocation given by the formula for jobs in $S_2$, if we replace $\Delta_{m}$ by $\Delta_{m}+\Delta_1\frac{o(1)}{o(m)}$ (the values $m_i$ and our assumption on the costs do not depend on the values $\Delta_i$). Increasing $\Delta_m$ only gives the adversary more options and does not change the cost of all jobs underaged. So by the induction assumption, every profile for jobs in $S_2$ where job $2$ starts at time $A_2$, has cost at most $\sum_{i\in S_2}u_ix_i=\sum_{i\in [k]}u_ix_i$. Now, consider an arbitrary profile $T$ for allocation $\{x_i\}_{i\in[k]}$ and recall that jobs $1$ to $m$ are all appointed at time $0$. If job $m$ is underaged in $T$, then setting $t_1=0$ reduces the completion times of jobs $1$ to $m$ by exactly $t_1$ and does not affect jobs after $m$. The resulting change in cost is $t_1(u_{m}-o(1)+o(m))\geq 0$, due to the assumption $u_{m_1}\geq o(1)-o(m_1)$. Since job $2$ starts at $A_2$ in the new profile, the cost of $T$ is at most $\sum_{i\in [k]}u_ix_i$. On the other hand, if job $m$ is overaged in $T$, let $d$ denote the delay with which job $m+1$ starts. Allocation $\{x_i\}_{i\in S_{m+1}}$ is given by the formula for jobs in $S_{m+1}$, so from the induction assumption and Lemma \ref{shift}, the total cost of jobs in $S_{m+1}$ is at most the cost of all jobs in $S_{m+1}$ taking maximum time with delay $d$. Each job $j<m$ has cost $o_jC_j$ and job $m$ is overaged, therefore both the cost of jobs $1$ to $m$ and the delay $d$ are non-decreasing in $t_1,\dots,t_m$. Increasing $t_1,\dots,t_m$ to their maximum values, we have that the cost of $T$ is at most the cost of all jobs overaged. This completes the induction.
	
	Now, to show that the allocation is an optimal solution to LP-2, consider the dual,
	\begin{eqnarray}
	&\max &\sum_{i=1}^{n} (\sum_{j=i}^n o(j) \Delta_j)p_i+\sum_{i=1}^n (\sum_{j=1}^i \Delta_j)r_i\nonumber\\
	&s.t. &(u_j+o(j)) \sum_{i\leq j} p_i + \sum_{i\geq j} r_i + s_j\leq u_j\quad \forall j\in[n]\label{dual}\\
	&& p_i\geq 0 \quad \forall i\in [n]\backslash \{1\};\quad -r_i,s_i\geq0\quad \forall i\in[n].\nonumber
	\end{eqnarray} 
	Consider the solution $p_1=\frac{u_{m_1}}{u_{m_1}+o({m_1})}$ and $p_i=\frac{u_{m_i}}{u_{m_i}+o({m_i})}-\frac{u_{{m_{i-1}}}}{u_{{m_{i-1}}}+o({m_{i-1}})}$ for every $i>1$ where $m_i\neq m_{i-1}$. Set all other variables to zero and note that when $m_i\neq m_{i-1}$ we have, $m_{i-1}=i-1$ and further, $\frac{u_{m_i}}{u_{m_i}+o({m_i})}-\frac{u_{{m_{i-1}}}}{u_{{m_{i-1}}}+o({m_{i-1}})}\geq0$. To check feasibility, it only remains to consider constraints (\ref{dual}). So fix arbitrary $k\in [n]$, and we have, $(u_k+o(k))\sum_{i\leq k}  p_i=(u_k+o(k))\frac{u_{m_k}}{u_{m_k}+o({m_k})}\leq u_k$, by definition of $m_k$. Finally, since the objective value of this dual solution is designed to  match the objective value of allocation $x_k=\frac{\sum_{i|k=m_i} \Delta_i o(i)}{u_k+o(k)} \quad \forall k\in[n]$ for the primal, we have a certificate of optimality for the primal.
	
\end{proof}
\section{Optimal Allocation for Convex Piece-wise Linear Costs}\label{ext2}
Given non-decreasing, convex piece-wise linear costs $o_i(\cdot),u_i(\cdot)$, we have the following lemmas.
\begin{lem}\label{shiftg}
	Given an arbitrary allocation $\{A_i\}_i$, where $t_i=\Delta_i$ for all $i$ is a worst-case and every job $i$ completes at or after time $A_{i+1}$ in this case (this holds, for instance, if $\sum_{j\leq i}x_j\leq \sum_{j\leq i}\Delta_j$ for all $i$). Recall that w.l.o.g., $A_1=0$ and consider a modified problem for the adversary, where the first job is always forced to start at a later time $t_0$ instead of time $A_1=0$. Then, $T=(t_i)$ is also a worst-case for the modified problem.
\end{lem}
\begin{proof}
	Let $T=(\Delta_i)$ denote the profile for all jobs taking maximum time. When job $1$ starts at time 0, denote the cost of job $i$ by $c(0,T,i)$. Since $T$ is a worst-case by assumption, we have $c(0,T)\geq c(0,Z)$ for every profile $Z$. Let $\hat{o}_i=\frac{c(t,T,i)-c(0,T,i)}{t}$. Then, for the modified setting where job $1$ starts at time $t>0$, we claim that $c(t,Z,i)\leq c(0,Z,i)+\hat{o}_it$ for every $Z$. To see this, suppose job 1 starts at 0 and consider two cases: (i) $i$ overaged in $Z$ and (ii) $i$ underaged in $Z$. In case (i), $i$ will still be overaged when job 1 starts at time $t$ and the completion time of $i$ can increase by at most $t$. In case (ii), if $i$ is still underaged when job 1 starts at time $t$, we are done. Else, $i$ becomes overaged. Clearly, the completion time of $i$ in $Z$ is no more than the completion time of $i$ in $T$. So due to the non-decreasing and convex nature of costs, the maximum increase in overage cost of $i$ in profile $Z$, due to late start $t$, is $\hat{o}_it$. 
	Therefore, $c(t,Z)\leq \sum_i (c(0,Z,i)+\hat{o}_i t)\leq \sum_i c(0,T,i)+\sum_i \hat{o}_i t= c(t,T)$.
\end{proof}


\begin{thm}
    Given non-decreasing, convex piece-wise linear costs $o_i(\cdot),u_i(\cdot)$ with $o_i(0)=u_i(0)=0$. Let $\olo_i ,\ulo_i$ denote the largest and smallest slope for $o_i(\cdot)$  and, $\olu_i,\ulu_i$ the largest and smallest slope for $u_i(\cdot)$. Then Algorithm \ref{optgalg} finds an optimal optimal allocation if the following is satisfied, \[\ulu_{i+1}\geq \olu_i, \text{ for all $i\in [n-1]$.}\] 
\end{thm}
\begin{proof} 
In case of ties, suppose the values $x_i$ are the minimum possible. This is w.l.o.g.\ since given the conditions on the cost functions, ties only occur if some costs are identically zero.

 Consider first the case of a single job i.e., $n=1$. Then clearly, $x_n=\arg\min_{x\geq 0} \max\{u_n(x),o_n(\Delta_n-x)\}$ is an optimal allocation. Further, since $u_n(0)=o_n(0)=0$, we have $u_n(x_n)=o_n(\Delta_n-x_n)$ and $x_n\leq \Delta_n$. We claim that more generally, the allocation $\{x_1,\dots,x_n\}$ output by the algorithm is such that, $x_i\leq \Delta_i\, \all i\in[n]$, and the worst-cases are given by all jobs taking zero time and all taking maximum time. To see this, first note that $u_i(0)=o_i(0)+\hat{o}_i(\Delta_i)=0$ and therefore, the minimum value of $\max\{u_i(x),o_i(\Delta_i-x)+\hat{o_i}(x)\}$ is attained for some $x_i\leq \Delta_i$ where, $u_i(x_i)=o_i(\Delta_i-x_i)+\hat{o_i}(x_i)$. From this we have that, 
 \[\sum_{i\in[n]\backslash [k-1]} u_i(x_i)=\sum_{i\in[n]\backslash [k-1]} \big(o_i(\Delta_i-x_i)+\hat{o_i}(x_i)\big)=\sum_{i\in[n]\backslash [k-1]} o_i\Big(\sum_{j=k}^{i}(\Delta_j-x_j)\Big).\]
 Here the last equality follows from definition of $\hat{o}_i(\cdot)$. This shows that the two cases (all jobs taking zero time, and all taking maximum time) have the same cost. We show that these are worst-cases by considering the set of jobs $\{k,\dots,n\}$ and performing induction on $k$. Clearly, for $k=n$ we have the desired. Assume the property of worst-cases for some $1<k<n$, and consider the case of $k-1$. Then, for the allocation $\{x_{k-1},\dots,x_n\}$, suppose that job $k-1$ is overaged in some worst-case. Then using the induction assumption and applying Lemma \ref{shiftg}, we have that all jobs overaged is a worst-case. Similarly, suppose that job $k-1$ is underaged (and thus takes zero time) in some worst-case. Then the case of all jobs taking zero time is a worst-case by the induction assumption. 
 
 Having characterized the worst-cases for allocation $\{x_i\}$, we now proceed to show optimality. In particular, the worst-case cost of $\{x_i\}$ is $\sum_i u_i(x_i)$. Given an arbitrary allocation $\{y_i\}_i$ and $k\in[n+1]$, let
 \[L_k(y)=\sum_{i<k}u_i(y_i)+\sum_{i\geq k}o_i\Big(\big(\sum_{l=k}^{i}(\Delta_l-y_l)\big)^+\Big).\]
 Consider the profile where jobs $1$ to $k-1$ take zero time and jobs $k$ to $n$ take maximum time. Each job $i<k$ is underaged with cost $u_i(y_i)$, job $k$ starts at $A_k$, and $C_i-A_{i+1}\geq \sum_{l=k}^{i}(\Delta_l-y_l)$ for every $i\geq k$. Since the costs are non-decreasing and $o_i(0)=0$, the cost of this profile is at least $L_k(y)$. Therefore, the worst-case cost of $\{y_i\}$ is at least $L_k(y)$ for every $k$. For the allocation $\{x_i\}$ output by the algorithm, recall that $x_i\leq\Delta_i$ and $\sum_{i\geq k}u_i(x_i)=\sum_{i\geq k}o_i(\sum_{l=k}^{i}(\Delta_l-x_l))$ for every $k$. Thus, $L_k(x)=\sum_iu_i(x_i)$ for every $k$.
 
 Now, fix an arbitrary allocation $\{y_i\}$, let $\delta_i=y_i-x_i$ and $\sigma_{j,i}=\sum_{l=j}^{i}\delta_l$. We call an index $k\in[n+1]$ \emph{good} if $\sigma_{k,i}\leq 0$ for all $i\geq k$ (note that $n+1$ is good), and let $k$ denote the smallest good index. We claim that $\sigma_{l,k-1}\geq 0$ for every $l<k$. Suppose not, and let $j$ be the largest $l<k$ with $\sigma_{l,k-1}<0$. Since $j$ is not good, there exists $i\geq j$ with $\sigma_{j,i}>0$. Note that $i\geq k$ is not possible since then $\sigma_{j,i}=\sigma_{j,k-1}+\sigma_{k,i}<0$. So we have $i< k-1$ and $\sigma_{i+1,k-1}=\sigma_{j,k-1}-\sigma_{j,i}<0$, contradicting the choice of $j$. Next, for every $l<k$ we have from convexity, $u_l(y_l)-u_l(x_l)\geq c_l\delta_l$ where $c_l=\ulu_l$ if $\delta_l\geq0$ and $c_l=\olu_l$ otherwise. Due to the condition $\ulu_{l+1}\geq\olu_l$, we have $0\leq c_1\leq c_2\leq\dots\leq c_{k-1}$. Therefore, using the claim above and rearranging the sum,
 \[\sum_{l<k}\big(u_l(y_l)-u_l(x_l)\big)\geq c_1\sigma_{1,k-1}+\sum_{l=2}^{k-1}(c_l-c_{l-1})\sigma_{l,k-1}\geq0.\]
 Finally, for every $i\geq k$ we have $\sum_{l=k}^{i}(\Delta_l-y_l)=\sum_{l=k}^{i}(\Delta_l-x_l)-\sigma_{k,i}\geq \sum_{l=k}^{i}(\Delta_l-x_l)\geq 0$. Combining everything, the worst-case cost of $\{y_i\}$ is at least $L_k(y)\geq L_k(x)=\sum_iu_i(x_i)$. Therefore, $x_i$ is an optimal allocation. 
\end{proof}
\section{Missing Proofs from Section \ref{rasog}}\label{apxrasog}

\begin{lem}
	Given an $\alpha$ approximation for $1||\sum f_j(C_j)$,  we have a $4\alpha n$ approximation for RASO. 
\end{lem}
\begin{proof}
	We claim that the value of (\ref{genshed}) is at most a factor $4n$ times the value of (\ref{problemn}). We show this by proving there always exists an order such that, given this order the cost of the allocation $x_j=\Delta_j \frac{o(j)}{o(j)+u_j}$, is within a factor $4n$ of (\ref{problemn}). Then the guarantee follows due to an $\alpha$ approximation for $1||\sum_jf_j(C_j)$.
	
	For jobs $j\in[n]$, let $y_j$ be an optimal allocation in an optimal order, such that all jobs underaged and all jobs overaged are worst-cases (see Lemma \ref{l1}). Recall that all jobs underaged and overaged are worst-cases for $y_j$ therefore, $\sum_j u_jy_j$ equals \eqref{problemn}. Further, let $C_j$ denote the completion time when all jobs are overaged and $B_j$ the appointed end time, for every job $j\in[n]$. Observe that for every job $j$ the maximum of the cost $u_jy_j$ (of underaging) or the cost $o_j(C_j-B_j)$ (of overaging), is at least $\Delta_j\frac{o_ju_j}{o_j+u_j}$ (optimal cost for an individual job). Further, $\Delta_j\frac{o_ju_j}{o_j+u_j}\geq \Delta_j\frac{\min\{o_j,u_j\}}{2}$. With $c_j=\max\{u_jy_j, o_j(C_j-B_j)\}$, we have, $\Delta_j\frac{o_ju_j}{o_j+u_j}\leq c_j\leq u_j y_j +o_j(C_j-B_j)$ and therefore, 
	\begin{equation}
	\Delta_j\frac{\min\{o_j,u_j\}}{2}\leq c_j \text{ and } \sum_j c_j\leq 2\sum_j u_j y_j.\label{cineq}    
	\end{equation}
	Now, let $S$ denote the set of jobs such that $c_j<\Delta_ju_j/2$. Then, for every job $j\in S$ we have, $c_j\geq \Delta_j o_j/2$. 
	
	Consider the order $\pi$ that schedules jobs not in $S$ in arbitrary order followed by jobs in $S$ in descending order of $o_j$. 
	Let $o'(j)=\sum_{k|\pi(k)\geq \pi(j)} o_k$. Then the cost of the formula allocation $x_j=\Delta_j \frac{o'(j)}{u_j+o'(j)}$ for this order is at most $4n$ times the optimal cost. To see this, recall that all jobs underaged is a worst-case for the formula allocation (Lemma \ref{formulaworstcase}) and observe that, 
	\begin{eqnarray*}
	\sum_{j\in S} u_jx_j &&\leq \sum_{j\in S} \Delta_j \min\{u_j,o'(j)\}\leq \sum_{j\in S} \Delta_j \min\{u_j,no_j\}\\
	&&\leq 2n \sum_{j\in S} c_j.
	\end{eqnarray*}
	Here the first inequality follows from $\frac{u_j o'(j)}{u_j+o'(j)}\leq\min\{u_j,o'(j)\}$, and the second inequality from the fact that jobs in $S$ are scheduled in decreasing order of overage costs. The last inequality follows from \eqref{cineq}. Finally, note that for a job $i\not\in S$, we have $u_ix_i\leq \Delta_iu_i\leq 2c_i$, regardless of the order. Thus, using \eqref{cineq} again, $\sum_j u_j x_j\leq 2n\sum_j c_j\leq 4n\sum_j u_j y_j$.
\end{proof}
The following example shows that there is a gap of $\Omega(n)$ between the values of (\ref{problemn}) and (\ref{genshed}). Thus, the analysis above is tight up to constant factors.

Consider jobs $1,\dots,2n$, and an order such that job $i$ is the $2n-i+1^{th}$ appointment scheduled. Let even indexed jobs have large underage cost $u_{2i}=M$ such that $\frac{M}{M+o(2n)}\to 1$, and small overage cost $o_{2i}=1\ll M$, for every $i\in[n]$. For the odd indexed jobs, let $o_{2i-1}=u_{2i-1}=1$ for every $i\in[n]$. Therefore, $\frac{u_{2i+1}}{u_{2i+1}+o(2i+1)}\leq\frac{u_{2i-1}}{u_{2i-1}+o(2i-1)}$ for every $i\in[n]$. Finally, let $\Delta_{2i}=1$ and $\Delta_{2i-1}=0$ for every $i\in[n]$. Now, for any two consecutive jobs $2i,2i-1$, using Theorem \ref{sidedspec} we have in the optimal allocation, $x_{2i}=0$ and $x_{2i-1}=(\Delta_{2i}o(2i)+\Delta_{2i-1}o(2i-1))\frac{1}{u_{2i-1}+o(2i-1)}$. 
Further, since $u_{2i-1}=o_{2i}=1$, we have $u_{2i-1}x_{2i-1}=\Delta_{2i}\frac{o(2i)u_{2i-1}}{u_{2i-1}+o(2i-1)}=1$. Thus, the cost of this order is $\sum_{i\in [2n]} u_i x_i=n$. This establishes an upper bound on the optimal value of (\ref{problemn}) for this instance. Next, consider an arbitrary order $\pi$ where we use the formula allocation and let $o'(i)=\sum_{k|\pi(k)\geq \pi(i)} o_k$. For this allocation, the cost of underaging an even indexed job is given by $\Delta_{2i}\frac{o'(2i)u_{2i}}{o'(2i)+u_{2i}}\to \Delta_{2i}o'(2i)$. Therefore, the total cost for this allocation is $\sum_{i\in[n]} \Delta_{2i}o'(2i)\geq \sum_{i\in[n]}i=n(n+1)/2$. Since this holds for every order, we have a multiplicative gap of $\Omega(n)$ between (\ref{problemn}) and (\ref{genshed}). 

 Using the same example, we can also show that no heuristic that orders jobs based on a simple ratio can beat $\Omega(n)$. 
 To see the claim, notice that every job in the above example is of one of two types. Denote the set of even indexed jobs as $E$ and $O$ is the set of odd indexed jobs. Therefore, any list order either schedules all jobs in $E$ first, followed by odd indexed jobs, or vice versa. For both orders, we can show that the optimal allocation has cost $\Omega(n^2)$. First, consider the order where type $O$ jobs are scheduled first. In this case, the optimal allocation is obtained by allotting zero time for all type $O$ jobs and using the formula allocation $x_i=\Delta_i\frac{o(i)}{o(i)+u_i}$ for type $E$ jobs. For this allocation, underaging all jobs is a worst-case with cost $\Theta( \sum_{\text{type }E}o(i))=\Theta(\sum_{i\in[n]}i)=\Theta(n^2)$. On the other hand, if the list order schedules all type $E$ jobs first, then consider the 2 approximate solution given by LP-1 for this order. In this allocation, all type $E$ jobs are allotted zero time and the first type $O$ job takes $n$ units of time. All other type $O$ jobs are allotted zero time. Overaging all type $E$ jobs, we have a cost of $O(n^2)$. Therefore, the worst-case cost for the LP allocation is $\Omega(n^2)$. Since this is guaranteed to be within a constant factor of the optimal, we get the desired.
 
\subsection{Possible Approach for Approximating RASO} \label{apxrasog2}

For fixed job order, recall the closed form optimal solution to LP-1: For every $i$, let $y_i=0$ if $u_i-o(i)\geq 0$ or $i\neq m_i$. For the remaining jobs, define $m_i=\arg\min_{j\geq i}u_j-o(j)$ and $y_i=\sum_{j|m_j=i}\Delta_j$. Let $M$ denote the set of all jobs $i$ such that $m_i=i$. Similar to the reduction in Theorem \ref{hard}, these jobs will act as a set of parallel ``machines''. Also, for every $i\not\in M$ let $o'(i)=u_{m_i}+\sum_{\{j \mid m_j=m_i, j \geq i, j\not\in M\}}o_j$. We will find that $o'(i)$ represents the completion time of job $i$ of machine $m_i$. Now, the worst-case cost of the LP optimal allocation is upper bounded by,
\begin{eqnarray}
 \sum_i [u_i y_i+o(i)(\Delta_i-y_i)]&&=\sum_{i\in M}u_i \sum_{j \mid m_j=i} \Delta_j+ \sum_{j\not\in M} o_j\sum_{\{k \mid m_k=m_j, k \leq j\}}\Delta_k\nonumber\\
 &&=\sum_{i \in M} u_i \sum_{j \mid m_j=i} \Delta_j+ \sum_{j\not\in M} \Delta_j (o'(j)-u_{m_j})\nonumber\\
 &&=\sum_{j \in M}u_j\Delta_j +\sum_{j \not\in M} o'(j)\Delta_j \label{ovcost}.
\end{eqnarray}
Thus, we have decomposition of this cost into two pieces; (i) A linear piece that represents fixed costs $u_i\Delta_i$ for each ``machine'' $i\in M$. (ii) The sum of weighted ``completion'' times $o'(i)$ (with weights $\Delta_i$) for regular (non-machine) jobs $i\not\in M$. Clearly, finding an order and allocation that minimizes \eqref{ovcost} to within a multiplicative factor of $\alpha$ would give us a $2\alpha$-approximation to RASO. One possible approach would be to break the problem down into two problems: (i) Finding a set of ``machine'' jobs and (ii) Finding a way to assign regular jobs to the machine jobs. A constant factor algorithm for both stages would yield a constant factor algorithm for the overall problem. And indeed, the latter problem of assigning regular jobs to machine jobs resembles a parallel machine scheduling problem (though the machines are not entirely identical), and one could use constant factor heuristics (such as the Smith's rule based heuristic in \cite{guchi,smithp}) to approximate the optimal value of (ii).  Unfortunately, it is unclear (to us) whether one can obtain a good approximation for (i) (the problem of finding the set of machine jobs) and therefore, more generally for minimizing \eqref{ovcost}.

 \section{Missing Proofs from Section \ref{rasohomo}}\label{apxrasohomo}
 \textbf{Sub-optimality of Algorithm \ref{localalg}:} Consider 3 jobs, with $\Delta_1=o_1=1$, $\Delta_2=15-\epsilon, \Delta_3=14$ and $o_2=o_3=5$. Let $\epsilon>0$ be an infinitesimally small quantity. Then, the order $(3,2,1)$ is locally optimal since, $\frac{\Delta_1}{o_1(o_1+1)}=\frac{1}{2}$ and $ \frac{\Delta_2}{o_2(o_2+1)}=\frac{1}{2}-O(\epsilon)$ and $\frac{\Delta_3}{o_3(o_3+1)}<\frac{1}{2}$. Further, $\frac{\Delta_2}{o_2(o_2+o_1+1)}> \frac{\Delta_3}{o_3(o_3+o_1+1)}$. However, notice that one can swap 1 and 2 without changing the value of the solution by more than $O(\epsilon)$, resulting in the order $(3,1,2)$ with nearly equal cost. However, this order is not locally optimal due to jobs 3 and 1, and by swapping them we have an order strictly better than the initial one.
 
 \begin{proof}[Proof of Lemma \ref{pairorder}]
	Consider the schedule $\pi'$ which swaps $i$ and $j$. Let $C_i,C'_i$ denote the completion time of job $i$ in $\pi$ and $\pi'$ respectively. Similarly, $C_j,C'_j$ for job $j$ and note that the completion times of all jobs except $i,j$ are unchanged. The difference in costs between $\pi'$ and $\pi$ is given by, 
	\begin{eqnarray*} 				
		\sum_k \Delta_k (f(C'_k)-f(C_k))=&&\frac{\Delta_i(o_i+o_x)}{o_i+o_x+1}+\frac{\Delta_j(o_j+o_i+o_x)}{o_j+o_i+o_x+1}\\
		&&-\frac{\Delta_i(o_i+o_j+o_x)}{o_i+o_j+o_x+1}-\frac{\Delta_j(o_j+o_x)}{o_j+o_x+1}\\
		=&&\Delta_j\frac{o_i}{(o_i+o_j+o_x+1)(o_j+o_x+1)}\\
		&&-\Delta_i\frac{o_j}{(o_i+o_j+o_x+1)(o_i+o_x+1)}\\
		=&& \frac{o_jo_i}{o_j+o_i+o_x+1}\Big(\frac{\Delta_j}{o_j(o_j+o_x+1)}-\frac{\Delta_i}{o_i(o_i+o_x+1)}\Big)
	\end{eqnarray*}
	Therefore, $\sum_k \Delta_k (f(C'_k)-f(C_k))\geq 0$ iff $\frac{\Delta_j}{o_j(o_j+o_x+1)}\geq\frac{\Delta_i}{o_i(o_i+o_x+1)}$.
\end{proof}

\begin{proof}[Proof of Lemma \ref{orderimpl}]
	Suppose $o_i> o_j$. We have for $x>0$,
	\begin{eqnarray*}
		\frac{\Delta_i}{o_i(o_i+1)}&&\geq\frac{\Delta_j}{o_j(o_j+1)}\\
		\frac{\Delta_i}{o_i(o_i+x+1)}\frac{o_i+x+1}{o_i+1}&&\geq\frac{\Delta_j}{o_j(o_j+x+1)}\frac{o_j+x+1}{o_j+1}\\
		\frac{\Delta_i}{o_i(o_i+x+1)}\Big(1+\frac{x}{o_i+1}\Big)&&\geq\frac{\Delta_j}{o_j(o_j+x+1)}\Big(1+\frac{x}{o_j+1}\Big)\\
		\implies	\frac{\Delta_i}{o_i(o_i+x+1)}&&>\frac{\Delta_j}{o_j(o_j+x+1)}\quad \text{(Contradiction)}.
	\end{eqnarray*}
	Therefore, $o_j\geq o_i$ as desired. Now starting with $\frac{\Delta_j}{o_j(o_j+x+1)}\geq\frac{\Delta_i}{o_i(o_i+x+1)}$ and repeating the steps above with the assumption $o_j\geq o_i$, we have that $\frac{\Delta_j}{o_j(o_j+x'+1)}\geq\frac{\Delta_i}{o_i(o_i+x'+1)}$ for every $x'\geq x$.
\end{proof}

\section{Analysis of C-Smith}\label{apxcsmith}

Consider list ordering heuristics that order jobs using a real valued function evaluated individually for each job, using only the job parameters ($p_j,w_j$). We call such a heuristic \emph{weight proportional} if for every job $i$, the function value is linearly proportional to the weight $w_i$. Reminiscent of Smith's rule, we call the function value for each job, the `ratio' for the job.
\begin{lem}\label{normalize}
	Given a weight proportional list order, any worst-case instance $I$ for the heuristic can be modified into another worst-case instance $I'$ where all jobs have the same ratio. Further, we can also assume that jobs have ratio 1. 
\end{lem}
\begin{proof}
	Proof simply generalizes Lemma 3.5 in \cite{smith} (and Lemma 21 in \cite{eptas}). Consider a list order  heuristic and suppose there exists a partition of the set of jobs into $k$ classes based on the ratio. All jobs within a class have the same ratio. Order the classes in decreasing order of the  representative ratio and note that the algorithm schedules jobs in class $j$ before jobs in class $j+1$ for every $1\leq j\leq k-1$. Let $A(j,I)$ denote the total cost of jobs in class $j$, given the ratio based order for instance $I$. Similarly, $OPT(j,I)$ denotes the total cost over jobs in the same class $j$ in the optimal order (in case of multiple optima, fix one for this proof). We use the shorthand $A(I),OPT(I)$ to denote the total costs over all jobs. Let $\alpha$ be such that in a worst-case instance $I$, we have, $A(I)= \alpha OPT(I)$. Given a worst-case instance $I$, define $i_I$ to be the class with the smallest index such that, $A(i,I)\geq \alpha OPT(i,I)$. Note that such an $i$ always exists, as otherwise, $\sum_{j=1}^k A(j,I)=A(I)< \alpha OPT(I)$. Now, let $I$ be a worst-case instance with the smallest number $d$, of distinct classes. We show that $d$ equals one by considering two cases and deriving a contradiction in each.
	
	First, if $i_I=1$, then we can simply consider a new instance $I'$ that consists of only jobs in class 1. We have that $OPT(I')\leq OPT(1,I)$ and $A(I')=A(1,I)$. Therefore, $A(I')\geq \alpha OPT(I')$ and we have contradiction to the minimality of $d$.
	
	Otherwise, assume $i_I>1$. Since the ordering heuristic is weight proportional, we can multiply the weights of all jobs in class $i_I$ by a factor $\delta >1$ such that the ratio for jobs in class $i_I$ now equals the ratio for jobs in class $i_I-1$. Call this instance $I_1$. We have $A(i,I_1)=\delta A(i,I)$ and therefore $A(I_1)=A(I)+(\delta-1)A(i,I)$. Similarly, $OPT(I_1)\leq OPT(I)+(\delta-1) OPT(i,I)$. Further, since $A(i,I)\geq \alpha OPT(i,I)$, we get, $A(I_1)\geq A(I) +(\delta-1)\alpha OPT(i,I)$. Finally, $A(I_1)\geq \alpha (OPT(I)+(\delta-1)OPT(i,I))\geq\alpha OPT(I_1)$. Therefore $I_1$ is also a worst-case instance, but with $d-1$ distinct classes, contradiction.  
	Finally, given an instance where all jobs have the same ratio, we can easily multiply all jobs weights by a scalar to normalize the ratio to one.
\end{proof}
\emph{Remark:} One might ask at this point if tie-breaking rules make a difference to the performance guarantee. Since we can always perturb weights adversarially by tiny amounts so that there are no ties but the costs are essentially unchanged, tie-breaking doesn't make a difference in theoretical performance.
\begin{lem}\label{order}
Given a set of jobs with identical C-Smith ratio, the worst-case order is to schedule jobs in ascending order of $p_i$ and the optimal order is to schedule in descending order of $p_i$.
\end{lem}
\begin{proof}
	Consider two consecutive jobs $k$ and $k+1$ in an optimal schedule. From Lemma \ref{pairorder}, we have $\frac{w_k}{p_k(p_k+C_{k-1}+1)}\geq \frac{w_{k+1}}{p_{k+1}(p_{k+1}+C_{k-1}+1)}$. Then applying Corollary \ref{ordercoro}, we have that $p_k\geq p_{k+1}$. This more generally implies jobs must be ordered in descending order of $p_k$ in the optimal schedule.
	
	Now consider the worst order and observe that if there is pair of consecutive jobs $k,k+1$ in this order such that $p_{k}>p_{k+1}$, then using Corollary \ref{ordercoro}, we have that swapping $k$ and $k+1$ would be strictly worse. Therefore, in the worst schedule jobs are ordered in ascending order of $p_i$. 
\end{proof}

Interestingly, in case of Smith's rule we have the opposite of the above lemma which we state as a fact below, since it will be useful later in the analysis.

\begin{fact}{[Lemma~3.6 in \cite{smith}, Proposition 20 in \cite{eptas}]} \label{fact}
Given a set of jobs with identical Smith's ratio, the worst-case order is to schedule jobs in descending order of $p_i$ and the optimal order is to schedule in ascending order of $p_i$.
\end{fact}

So far, our simplification implies that we can assume all jobs have C-Smith ratio of 1. Given this, we also know the worst possible (C-Smith) order and the optimal order. Now, in a worst-case instance we could have infinitely many jobs. In particular, we could have infinitely many infinitesimally small jobs. All such jobs that have the same processing time can be assumed to be scheduled contiguously in both C-Smith and the optimal order. Therefore, infinitely many of such infinitesimally short jobs can be amalgamated into a special job with finite length. We call such a job, a \emph{line} job. The notion was introduced in \cite{eptas,smith}, for the analysis of Smith's rule.

\emph{A line job of length $q$ consists of infinitely many infinitesimally small jobs, each with weight and processing time equal to $\epsilon$ (so the Smith's ratio is 1). The total number of jobs is $q/\epsilon$ with $\epsilon\to 0$. The total weight of a line job is also $q$, same as its processing time.}

We use the same notion here, but we let the weight of each small job be $\epsilon(\epsilon+1)$ so that its C-Smith ratio is 1. However, the $\epsilon^2$ term in the weight disappears upon integration and the cumulative weight of a line job of length $q$ in our setting is also $q$. For the rest of the analysis, we consider a worst-case given by the following proposition.
\begin{proposition}
For some integer $n\geq 0$, there exists a worst-case instance for C-Smith with $n+1$ jobs, $n$ of which are regular jobs that have processing times, $p_1\geq p_2\geq \dots \geq p_n$, and the $n+1$ th job is a line job of length $q$. All jobs have C-Smith ratio 1. Therefore, regular job $i$ has weight $p_i(p_i+1)$ and the line job has weight $q$. 
\end{proposition}

Recall, our concave function over completion time -- $f(C_j)=\frac{C_j}{1+C_j}$. The cost incurred due to a line job that runs from time $t$ to $t+p$ is given by, 
\[ F(t,t+p):=\int_{t}^{t+p} f(x)dx=\int_{t}^{t+p} \Big(1-\frac{1}{x+1}\Big)dx = \quad p-\ln\Big(\frac{p}{t+1}+1\Big) .\]
Observe that $F(t,t+p)+F(t+p,t+p+q)=F(t,t+p+q)$ for every $t,p,q>0$, and this telescoping property of $F(.)$ will be useful later. Further, since $f(.)$ is monotonically increasing we have,
\begin{eqnarray}
p\frac{t}{t+1}\leq F(t,t+p)\leq p\frac{t+p}{t+p+1}\label{lineineq}.
\end{eqnarray}
Let $p(i)=\sum_{k\leq i} p_k$ and $q(i)=q+ \sum_{k\geq i}p_k$ with $p(0)=0$ and $q(n+1)=q$ (note the difference between the sums in $p(i)$ and $q(i)$). Note that weight $w_i$ of job $i$ is given by $p_i(p_i+1)$. 
Now, the optimal order schedules jobs 1 through $n$ in that order and ends with the line job. The worst-case schedule is the exact reverse, beginning with the line job, followed by job $n$ and ending with job 1. We let $p_1>0$, otherwise the optimal and C-Smith order are identical for all values of $q$ and we have a ratio of 1. Given these definitions, the worst-case ratio $W(n)$, for instances with $n$ regular jobs and a line job, all with C-Smith ratio unity, can be written as,
\begin{eqnarray}
N(q,p_1,\dots,p_n)&&=F(0,q)+\sum_{i=1}^n p_i(p_i+1)\frac{q(i)}{q(i)+1}\label{N}\\
D(q,p_1,\dots,p_n)&&=F(p(n),p(n)+q)+\sum_{i=1}^n p_i(p_i+1)\frac{p(i)}{p(i)+1}\label{D}\\
&&W(n)=\max_{q,p_1,\dots,p_n\geq 0} \frac{N}{D}(q,p_1,\dots,p_n)\nonumber.
\end{eqnarray}
Clearly $W(n)\leq W(n+1)$ for all $n$. Interestingly, for Smith's rule one could simplify the worst-case instance even further. In fact, there exists a worst-case which includes just one regular job and a line job. More formally, if we use $S(n)$ to denote the worst-case ratio for Smith's rule with $n$ regular jobs, \cite{smith,eptas} showed that $S(n)=S(1)$ for all $n\geq 1$. The worst-case for C-Smith is more involved and in particular, it is not clear if there is a worst-case with a small number of regular jobs.
In fact, we show that $W(1)<W(2)<W(3)<W(4)$, implying that every worst-case has at least $4$ regular jobs. Nonetheless, we will show that $ 1.06036<W(4)$ and $W(n)<1.06043$ for every $n\geq 2$, effectively showing $W(4)$ is reasonably close to the true approximation guarantee. To establish the upper bound on $W(n)$, we introduce a family of non-linear optimization problems where we can control the number of variables. 

In order to construct this family of upper bounds, we establish and exploit constraints on the \emph{worst-case} value of variables $q,p_1,\dots,p_n$ i.e., values that attain objective value $W(n)$. Further, the same bounds and constraints also greatly help when using Couenne to actually solve instances of the non-convex problems in reasonable time. To give a quick overview of the rest of the analysis: First, we show a lower bound on $p_1$ that allows us to find $W(1)$. Without the lower bound, the solver does not converge due to numerical issues related to dividing two very small quantities as $p_1,q\to 0$. Then, we show that $W(2)>W(1)$ by finding a lower bound on $W(2)$ that is larger than $W(1)$. This implies that there are at least two regular jobs in every worst-case. Given this observation, we then establish bounds and constraints on other variables. Finally, we combine these conditions to establish our family of upper bounds, which we then solve numerically. 


\begin{lem}\label{lbp}
	For every $n$, in a worst-case instance with $n$ regular jobs (i.e., one that attains $W(n)$ defined above), we have $p_1\geq W(n)-1$. In particular, $p_1\geq \beta'-1$ for every $\beta'\leq W(n)$.
\end{lem}
\begin{proof}
	Let us assume $p_1< W(n)-1$ and derive a contradiction.
	Considering the numerator (\ref{N}) for $n$ jobs we have,$$q-\ln(q+1) +\sum_{i=1}^n p_i(p_i+1)\frac{q(i)}{q(i)+1}< W(n) \Big(q-\ln(q+1) +\sum_{i=1}^n p_i\frac{q(i)}{q(i)+1}\Big):=W(n) N' ,$$
	where the inequality follows from $p_i+1\leq p_1+1< W(n)$. On the other hand, the denominator (\ref{D}), $$q-\ln\Big(\frac{q}{p(n)+1}+1\Big) +\sum_{i=1}^n p_i(p_i+1)\frac{p(i)}{p(i)+1}\geq q-\ln\Big(\frac{q}{p(n)+1}+1\Big) +\sum_{i=1}^n p_i\frac{p(i)}{p(i)+1}:=D'.$$
	The key now is to observe that $\frac{N'}{D'}$ can be interpreted as the ratio of the optimal order divided by the worst-case order for Smith's rule, on an instance where all jobs have Smith's ratio 1 i.e., job $i$ takes time $p_i$ and has weight $p_i$ for all regular jobs $i\in[n]$ and in addition we have a line job of length and weight $q$. From Fact \ref{fact}, $N'$ represents the optimal cost and $D'$ the worst-case cost. Therefore, $\frac{N'}{D'}\leq 1$. This gives a contradiction since, $\frac{N}{D}<W(n) \frac{N'}{D'}\leq W(n)$.
\end{proof}

Now that we have a lower bound on $p_1$, we evaluate $W(1)$ using Couenne (with $\beta'=1.05$, which is a valid lower bound on $W(1)$) and find,
\begin{eqnarray}
W(1)=\max_{p\geq \beta'-1, \, q\geq 0} 
\frac{q-\ln(q +1)+p(p+1)\frac{p+q}{p+q+1}}{q-\ln(\frac{q}{p+1}+1)+p^2}= 1.0567 \text{ when } p=2.8943, q=8.0775. \label{betalb}
\end{eqnarray} 
Further, we also have $W(2)>1.0603>W(1)$ for $q=7.7309, p_1=2.3001,p_2=1.8336$. 
We now establish more constraints on the variables but first, consider the following simple and useful observation.

\begin{lem}\label{ratio}
	For any value $\delta\geq 1$, given $\frac{N}{D}=\delta$, and $n,d\geq 0$, we have the following:
	\begin{enumerate}
		\item If $n\geq \delta d$ then $\frac{N+n}{D+d}\geq \delta$.
		\item If $n\leq \delta d$ and $D-d>0$ then $\frac{N-n}{D-d}\geq \delta$.
	\end{enumerate}
	In both cases, we have equality iff $n=\delta d$.
\end{lem}
\begin{proof} \textit{1.} follows from $\frac{N+n}{D+d}\geq \frac{N+\delta d}{D+d}=\delta$ and \textit{2.} from $\frac{N-n}{D-d}\geq \frac{N-\delta d}{D-d}=\delta$. 
\end{proof}
We now show a lower bound on $q$, which is useful in proving several subsequent results.
\begin{lem}\label{lbq}
	For $n\geq 2$, in a worst-case with $p_n>0$ we have, $q> p(n-1)$.
\end{lem}
\begin{proof}
	Consider the order given by C-Smith and the optimal order and suppose we turn job $n$ (for $n\geq 2$) into a line job with length and weight $p_n$, keeping everything else as is in both orders. The change in cost (new cost - old cost) for the C-Smith schedule is,
	\begin{eqnarray}
	\delta_N=-p_n(p_n+1)\frac{q(n)}{q(n) +1}+F(q(n+1),q(n)).
	\end{eqnarray}
	Similarly, the change in cost of the optimal order is,
	\begin{eqnarray}
	\delta_D=-p_n(p_n+1)\frac{p(n)}{p(n) +1}+F(p(n-1),p(n)).
	\end{eqnarray}
	Using (\ref{lineineq}), it is easy to see that both $\delta_N$ and $\delta_D$ are negative. 
	So if $-\delta_N\leq -\delta_D$, then we have a strictly worse instance due to Lemma \ref{ratio} \textit{(2.)}. 
	
	We claim $-\delta_N\leq -\delta_D$ as long as $q(n)\leq p(n)$ or $q\leq p(n-1)$. Therefore, we can assume that $q> p(n-1)$ in the worst-case. To see the claim, treat $\delta_N$ as a function of $q$ i.e., $\delta_N(q)$. Clearly, 
	$\delta_N(p(n-1))=\delta_D$. Moreover, $\frac{d\delta_N}{dq}=-p_n(p_n+1)\frac{1}{(q+p_n+1)^2}+\frac{p_n}{(q+1)(q+p_n+1)}=\frac{-qp^2_n}{(q+1)(q+p_n+1)^2}\leq 0$, therefore $\delta_N(q)\geq \delta_D$ for $0\leq q\leq p(n-1)$ and $\delta_N(q)\leq \delta_D$ for $q>p(n-1)$.
\end{proof}
The following is a useful property of $F$, which combined with the above lower bound on $q$, will help us establish a lemma that is key to our upper bounds on $W(n)$.
\begin{lem}\label{Fprop}
	Given $t_2,t_1,p\geq 0$, we have $p\frac{t_1+p}{t_1+p+1}-F(t_1,t_1+p)> p\frac{t_2+p}{t_2+p+1}-F(t_2,t_2+p)$ iff $t_2>t_1$.
\end{lem}
\begin{proof}
	Observe that $p\frac{t+p}{t+p+1}-F(t,t+p)=\int_{t}^{t+p} (\frac{t+p}{t+p+1}-\frac{t+x}{t+x+1})dx\geq 0$. Since $f$ is strictly concave, we have $t_2>t_1$ iff $\frac{t_1+p}{t_1+p+1}-\frac{t_1+x}{t_1+x+1}> \frac{t_2+p}{t_2+p+1}-\frac{t_2+x}{t_2+x+1}$ for every $x\geq0$.  
\end{proof}

\begin{lem}\label{amalgam}
Let $p(i)=\sum_{j=1}^i p_j$, and consider the ratio,
	\[r=\frac{\lambda p_n\frac{q+p_n}{q+p_n+1}+P}{\lambda p_n \frac{p(n)}{p(n)+1}+Q}.\] Then, given $r,\lambda \geq 1$, $P,Q,q,p_n\geq 0$, and $q+p_n>p(n)\geq p_n$ we have,
	\[\frac{\lambda F(q,q+p_n)+P}{\lambda F(p(n)-p_n,p(n))+Q}>r.\]
\end{lem}
\begin{proof}
 The claim is trivial when $p_n=0$, so let us assume $p_n>0$. 
	Suppose in the definition of $r$, we replace the term $p_n \frac{q+p_n}{q+p_n+1}$ in the numerator by $F(q,q+p_n)$ 
	and replace $p_n \frac{p(n)}{p(n)+1}$ in the denominator by $F(p(n-1),p(n))$. 
	Then, the numerator decreases (old cost - new cost) by $\delta_N=\lambda (p_n \frac{q+p_n}{q+p_n+1}-F(q,q+p_n))$. 
	The denominator decreases by $\delta_D=\lambda (p_n\frac{p(n)}{p(n)+1}-F(p(n-1),p(n)))$. Using (\ref{lineineq}), clearly $\delta_N,\delta_D\geq 0$. Moreover, since we are given $q>p(n-1)$, using Lemma \ref{Fprop} we have, $\delta_D>\delta_N$. 
	Since the original ratio is at least 1, we have a strictly larger ratio after this operation due to Lemma \ref{ratio} \textit{(2.)}. 
\end{proof}
Using the above lemma, we now show an upper bound on the total time taken by regular jobs, $p(n)$.
\begin{lem}\label{nlsumbound}
	Let $p(i)=\sum_{j=1}^i p_j$. Given a worst-case instance $q,p_1,\dots,p_n$ for $W(n)$ (or for $W(n,\beta')$ defined below, with $\beta'\leq W(n)$), with value $\delta=\frac{N}{D}(q,p_1,\dots,p_n)$, $n\geq 2$, $p_n>0$, $q>p(n-1)$, and $p_i\geq p_{i+1}$ for every $i\in[n-1]$, we have, \[\frac{(p_n+q)(p(n)+1)}{(p_n+q+1)p(n)}\geq \delta .\]
\end{lem}
\begin{proof}
	Consider the numerator $N$ and denominator $D$ from (\ref{N}),(\ref{D}). We assume that $\frac{\frac{p_n+q}{p_n+q+1}}{\frac{p(n)}{p(n)+1}}< \delta$ and derive a contradiction. Using Lemma \ref{ratio} \textit{(2.)}, we have that replacing $p_n(p_n+1)\frac{q(n)}{q(n)+1}$ in $N$ and $p_n(p_n+1)\frac{p(n)}{p(n)+1}$ in $D$ by $p_n\frac{q(n)}{q(n)+1}$ and $p_n\frac{p(n)}{p(n)+1}$ respectively, increases the objective. And finally, using Lemma \ref{amalgam} (with $\lambda=1$) we have that turning job $n$ into a line job of length $p_n$ increases the ratio, giving us a contradiction since the instance with line job of length $q+p_n$ and regular jobs $1,\dots,n-1$ has ratio larger than $\delta$, while $\delta\geq W(n-1)$ (this instance also satisfies the constraints of \eqref{nlp1} below). Note that we satisfy all conditions necessary to use Lemma \ref{amalgam} since the ratio $\delta$ is at least 1, $q>p(n-1)$ and other terms in the numerator and denominator are unaffected when converting $n$ into a line job. 
\end{proof}
\begin{coro}\label{sumbound}
	Given a worst-case solution with $n\geq 2$ regular jobs ($p_n>0$), $p(n)\leq \frac{1}{W(n)-1}$.
\end{coro}
\begin{proof}
	Observe that $\frac{(p_n+q)(p(n)+1)}{(p_n+q+1)p(n)}\leq \frac{p(n)+1}{p(n)}$. Then a direct application of Lemma \ref{nlsumbound} gives the desired.
\end{proof}

Now, suppose that for the case of $n$ regular jobs there exists a worst-case with $p_n>0$, then combining several of the results above we now have,
\begin{eqnarray}
&&W(n)=\max_{q,p_1,\dots,p_n} \frac{N}{D}(q,p_1,\dots,p_n)\label{nlp1}\\
&\text {s.t. } &q\geq p(n-1) \nonumber\\
&&p_{i}\geq p_{i+1}\quad \forall 1\leq i\leq n-1\nonumber\\
&& q, p_1\geq \beta-1\nonumber\\ 
&&p(n) \leq \frac{1}{\beta-1} \nonumber\\
&& q\geq 0; \quad p_i\geq 0\quad \forall i\nonumber
\end{eqnarray}
Due to Lemmas \ref{lbp}, \ref{lbq} and Corollary \ref{sumbound}, a worst-case instance with $p_n>0$ satisfies these constraints even if $\beta$ is replaced by any value $\beta'\leq W(n)$. The true approximation factor is given by $\lim_{n\to\infty} W(n)=\beta$. The limit exists since the sequence $W(n)$ is non-decreasing. However, solving (\ref{nlp1}) is non-trivial since the objective is provably non-convex. In fact, it can be shown that optimal solutions to $W(n-1)$ are saddle points for $W(n)$. Apart from the difficulty inherent in the structure of the optimization problem, another issue is that the formulation above uses the knowledge $\beta$ of the optimum value. However, the following proposition helps deal with this issue.
\begin{proposition}\label{propw}
For any $n$, given two values $\beta_1,\beta_2$ with $\beta_1\geq \beta_2\geq 1$, the optimum values $W(n,\beta_i)$, obtained by solving program \eqref{nlp1} with value $\beta_i$ in the linear constraints, obey $W(n,\beta_1)\leq W(n,\beta_2).$ Further, if the true optimum $\beta$ is lower bounded by $\beta_i$, then $\lim_{n\to\infty} W(n,\beta_i)=\beta$ and $W(n,\beta_i)\leq \beta$.
\end{proposition}
\begin{proof}
The proposition follows from the fact that smaller values of $\beta_i$ lead to a larger feasible region. As $n\to\infty$, the optimum value of the program converges to $\beta$, even absent the linear inequalities (therefore, $\lim_{n\to \infty} W(n,1)=\beta$). Finally, it is easily seen that $W(n,\beta_i)\leq W(n+1,\beta_i)$ i.e., $W(n,\beta_i)$ is non-decreasing in $n$ for any $\beta_i >1$.
\end{proof}
Now from \eqref{betalb}, we have a lower bound on $\beta$ in fact, we have a more refined lower bound from the fact that $W(2)>1.0603$. Solving the program using this lower bound, we find that $W(4,1.0603)>W(3,1.0603)>1.0603$, and the program above is already intractable to solve for $W(5,1.0603)$ and beyond using Couenne (solver does not converge within a week on a computing cluster). So while we have more refined lower bounds on the true guarantee $\beta$ (since $\beta \geq W(4,1.0603)$), in order to find upper bounds we need to account for the possibly large number of regular jobs in a worst-case instance.

\begin{lem} \label{ub}
	For $n\geq 2$, given a lower bound $\beta'>1$ on $\beta$, suppose that there is a worst-case instance for $W(n,\beta')$ with $p_n>0$. Then for every value of parameter $m$, such that $1\leq m\leq n-1$, $\beta$ is upper bounded by the optimal value of the following non-linear program,
	\small
	\begin{eqnarray}
	&\max  &u(q,r,p_1,\dots,p_m) \nonumber\\
	&&=\frac{F(0,q)+(p_{m}+1)F(q,q+r)+\sum_{i=1}^{m}p_i(p_i+1)\frac{q'(i)}{q'(i)+1}}{F(p(m)+r,p(m)+r+q)+(p_{m}+1)F(p(m),p(m)+r)+\sum_{i=1}^{m}p_i(p_i+1)\frac{p(i)}{p(i)+1}}\nonumber\\
	&\text {s.t. }&q\geq p(m)\label{nlp2} \\
	&&p_{i}\geq p_{i+1}\quad \forall 1\leq i\leq m-1\nonumber\\
	&& q, p_1\geq \beta'-1\nonumber\\ 
	&&p(m)+r \leq \frac{1}{\beta'-1} \nonumber\\
	&&q'(i)=q+r+\sum_{i\leq j\leq m}p_j\quad \forall i\leq m\nonumber\\
	&& q,r\geq 0; \quad p_i\geq 0\quad  \forall i\nonumber
	\end{eqnarray}
	\normalsize
\end{lem} 
\begin{proof}
	Let $U(m)$ denote the optimal value of the above program and note that the program is independent of $n$. So, consider any $n\geq 2$ such that there is an optimal solution $q,p_1,\dots,p_{n}$ for $W(n,\beta')$ with $p_{n}>0$.
	We show that for every $m\leq n-1$, $u(q,\sum_{i=m+1}^{n} p_i,p_1,\dots,p_m)$ upper bounds $W(n,\beta')$. Since the choice of $n$ is arbitrary (as long as $p_{n}>0$ and $n\geq 2$), and $U(m)$ does not depend on $n$, from Proposition \ref{propw} we have, $U(m)\geq \beta$ for every $m$. Note, the constraints in the program (\ref{nlp2}) follow almost directly from their counterparts in (\ref{nlp1}).   
	
	Fix arbitrary $m\geq 1$. To show $u(q,\sum_{i=m+1}^{n} p_i,p_1,\dots,p_m)\geq W(n,\beta')$ for $n\geq m+1$, let $p(i)=\sum_{j =1}^i p_j$ and $q(i)=q+\sum_{j=i}^n p_j$ as before. From Lemma \ref{nlsumbound}, we have for every $i\leq n$, \[\frac{\frac{q(i)}{q(i)+1}}{\frac{p(i)}{p(i)+1}}\geq \frac{\frac{q(n)}{q(n)+1}}{\frac{p(n)}{p(n)+1}}\geq W(n,\beta'). \] 
	So for every $i\geq m+1$, replacing $p_i(p_i+1)\frac{q(i)}{q(i)+1}$ in the numerator (\ref{N}) and $p_i(p_i+1)\frac{p(i)}{p(i)+1}$ in the denominator (\ref{D}) by $(p_{m}+1)p_i\frac{q(i)}{q(i)+1}$ and $(p_{m}+1)p_i\frac{p(i)}{p(i)+1}$ respectively, we get,
	\small
	$$V(q,p_1,\dots,p_n)=\frac{F(0,q)+(1+p_m)\sum_{i=m+1}^n p_i\frac{q(i)}{q(i)+1}+ \sum_{i=1}^m p_i(p_i+1)\frac{q(i)}{q(i)+1}}{F(p(n),p(n)+q)+(1+p_m)\sum_{i=m+1}^n p_i\frac{p(i)}{p(i)+1}+\sum_{i=1}^m p_i(p_i+1)\frac{p(i)}{p(i)+1}}.$$
	\normalsize
	Due to Lemma \ref{ratio} {\it (1.)}, we have $V(q,p_1,\dots,p_n)\geq W(n,\beta')$. We now repeatedly apply Lemma \ref{amalgam} on $V$ with $\lambda=1+p_{m}$ to obtain,
	\footnotesize
	$$\frac{F(0,q)+(1+p_m)\sum_{i=m+1}^n F(q(i+1),q(i))+ \sum_{i=1}^m p_i(p_i+1)\frac{q(i)}{q(i)+1}}{F(p(n),p(n)+q)+(1+p_m)\sum_{i=m+1}^n F(p(i-1),p(i))+\sum_{i=1}^m p_i(p_i+1)\frac{p(i)}{p(i)+1}}>V(q,p_1,\dots,p_n).$$
	\normalsize
	Here we applied the lemma once for every $i\geq m+1$. Note that, a crucial requirement for applying Lemma \ref{amalgam} requires $q+p_n>p(n)$, and this follows from the constraints (which are due to Lemma \ref{lbq}). Now, since the terms $F(q(i+1),q(i)))$ and $F(p(i-1),p(i))$ telescope, the LHS above is exactly $u(q,r,p_1,\dots,p_m)$ for $r=\sum_{i\geq m+1} p_i$. 
	To finish the proof note that  $(q,\sum_{i=m+1}^{n} p_i,p_1,\dots,p_m)$, is feasible for \eqref{nlp2} since given the substitution $r=\sum_{i\geq m+1} p_i$, the constraints in (\ref{nlp2}) follow directly from their counterparts in (\ref{nlp1}).
\end{proof}
We are now ready to combine everything together and prove Theorem \ref{maint}.
\begin{proof} \emph{(of Theorem \ref{maint})}
Solving program (\ref{nlp1}) using Couenne for $n=4$ and with $\beta'=1.0603$ in the constraints, we find $W(4)=W(4,\beta')=1.060369>W(3).$ 
	 Therefore, every worst-case has at least 4 regular jobs. We deal with the possibility of more than 4 jobs by using the family of upper bounds for $m=3$. Solving the 5 variable problem using Couenne, we get $U(3)\approx1.060429$ and therefore, $\beta<1.06043$.
\end{proof}
\end{appendix}
\end{document}